%% file: Abund-Surf-imperf05.tex
\documentclass[12pt,twoside]{amsart}
\usepackage{amssymb}
\usepackage{amscd}
\usepackage{color}
\usepackage{hyperref}

\usepackage[normalem]{ulem} 
\newcommand\redout{\bgroup\markoverwith
{\textcolor{red}{\rule[.5ex]{2pt}{0.4pt}}}\ULon}

\title[Abundance theorem]
{Abundance theorem for surfaces over imperfect fields} 
\author{Hiromu Tanaka} 
\address{Graduate School of Mathematical Sciences, 
The University of Tokyo, 
3-8-1 Komaba, Meguro-ku, Tokyo 153-8914, JAPAN} 
\email{tanaka@ms.u-tokyo.ac.jp}
\newcommand{\red}[0]{{\operatorname{red}}}

\newcommand{\Spec}[0]{{\operatorname{Spec}}}

\newcommand{\Supp}[0]{{\operatorname{Supp}}}

\newcommand{\Ex}[0]{{\operatorname{Ex}}}
\newtheorem{thm}{Theorem}[section]
\newtheorem{lem}[thm]{Lemma}

\newtheorem{prop}[thm]{Proposition}

\newtheorem{step}{Step}

\theoremstyle{definition}

\newtheorem{dfn}[thm]{Definition}

\newtheorem{rem}[thm]{Remark}
\newtheorem*{ack}{Acknowledgements}      
\newtheorem{nota}[thm]{Notation}

\newcommand{\MO}{\mathcal{O}}
\newcommand{\R}{\mathbb{R}}
\newcommand{\Q}{\mathbb{Q}}
\newcommand{\Z}{\mathbb{Z}}

\begin{document}

\maketitle

\begin{abstract}
In this paper, we show the abundance theorem for log canonical surfaces 
over fields of positive characteristic.
\end{abstract}

\tableofcontents

%%%%%%%%%%%%%%%

\input section1.tex

\input section2.tex

\input section3.tex

\input section4.tex

\input section5.tex

\end{document}

%% file: section1.tex
\section{Introduction}

The Italian school of algebraic geometry in the early twentieth century 
established a classification theory for smooth complex surfaces, 
which was later generalised by Kodaira, Shafarevich and Bombieri--Mumford. 
In particular, Bombieri--Mumford showed the abundance theorem 
for a smooth projective surface $X$ over an algebraically closed field $k$, 
i.e. if $K_X$ is nef, then it is semi-ample. 
After that, Fujita succeeded to extend it to a logarithmic case \cite{Fuj84}, 
which was extended to $\Q$-factorial surfaces by Fujino and the author 
(\cite{Fujino}, \cite{minimal}). 
As a consequence, Fujita established the abundance theorem 
for log canonical surfaces over any algebraically closed fields. 
The purpose of this paper is to prove 
the abundance theorem for log canonical surfaces 
over any field of positive characteristic. 
%which does not seem to be known even for regular surfaces. 

\begin{thm}\label{0abundance}
Let $k$ be a field of positive characteristic. 
Let $(X, \Delta)$ be a log canonical surface over $k$, where $\Delta$ is an effective $\R$-divisor. 
Let $f:X \to S$ be a proper $k$-morphism to a scheme $S$ which is separated and of finite type over $k$. 
If $K_X+\Delta$ is $f$-nef, then $K_X+\Delta$ is $f$-semi-ample. 
\end{thm}

If $k$ is a perfect field, then Theorem~\ref{0abundance} 
immediately follows from the case when $k$ is algebraically closed. 
Thus the essential difficulty appears only when $k$ is an imperfect field. 

There are some advantages to study varieties over imperfect fields, 
even if one works over an algebraically closed field of positive characteristic. 
For instance, given a fibration between smooth varieties of positive characteristic, general fibres are not smooth in general, 
whilst the generic fibre is always a regular scheme (e.g. quasi-elliptic surfaces).

On the other hand, 
there is a possibility to find new phenomena 
appearing only in positive characteristic 
by studying geometry over imperfect fields. 
Indeed, Schr\"{o}er succeeded to discover a three-dimensional del Pezzo fibration 
$f:X \to C$ with $R^1f_*\MO_X \neq 0$ by studying 
del Pezzo surfaces over imperfect fields (cf. \cite{Schroer}, \cite{Maddock}).

\subsection{Sketch of the proof}

Let us look at some of the ideas of the proof of Theorem~\ref{0abundance}. 
To this end, we assume that 
$S=\Spec\,k$ and $\Delta$ is a $\Q$-divisor throughout this subsection.

\subsubsection{Klt case}

If $k$ is perfect, then 
we can show Theorem~\ref{0abundance} just by taking the base change to the algebraic closure. 
However, if the base field is imperfect, the situation is subtler. 
First, the normality and reducedness might break up 
by taking the base change. 
Second, for the case when $k$ is algebraically closed, 
the proof of the abundance theorem for smooth surfaces 
depends on Noether's formula and Albanese morphisms. 
These techniques can not be used in general for varieties over non-closed fields. 
In particular it seems to be difficult 
to imitate the proof for the case over algebraically closed fields.

Let us overview the proof of Theorem~\ref{0abundance}. 
To clarify the idea, 
we treat the following typical situation: 
$X$ is a projective regular surface over a field of characteristic $p>0$ and 
$K_X \equiv 0$. 
What we want is to show that $K_X \sim_{\Q} 0$. 
First we can assume that $k$ is separably closed by taking the base change to the separable closure. 
Such a base change is harmless because, for example, a finite separable extension is \'etale. 
Second, we take the normalisation $Y$ of $(X\times_k \overline{k})_{\red}$. 
%where $(X\times_k \overline{k})_{\red}$ is the reduced structure of $X\times_k \overline{k}$. 
Set $f:Y \to X$ to be the induced morphism: 
$$f:Y \to (X\times_k \overline{k})_{\red} \to X\times_k \overline{k} \to X.$$
We can check that $Y$ is $\Q$-factorial (Lemma~\ref{purely-basic}(3)). 
By \cite[Theorem~1.1]{purely} (cf. Theorem~\ref{purely}), 
we can find an effective $\Z$-divisor $E$ on $Y$ such that 
$$K_Y+E=f^*K_X.$$
If $E$ is a prime divisor, then we can apply the following known result.

\begin{thm}[Theorem~1.2 of \cite{minimal}]\label{Q-fac-abundance}
Let $Y$ be a projective normal $\Q$-factorial surface 
over an algebraically closed field of characteristic $p>0$. 
Let $\Delta_Y$ be a $\Q$-divisor on $Y$ 
whose coefficients are contained in $[0, 1]$. 
If $K_Y+\Delta_Y$ is nef, then $K_Y+\Delta_Y$ is semi-ample. 
\end{thm}

Thus, we focus on 
the simplest but non-trivial case: $E=2C$ where $C$ is a prime divisor. 
Then we get    
$$K_Y+2C=f^*K_X.$$
%To apply Theorem~\ref{Q-fac-abundance}, it suffices to show that $C$ is semi-ample. 
Set $C_X:=f(C)$. 
We have the three cases: $C^2>0$, $C^2=0$, and $C^2<0$. 
If $C^2>0$, then $C$ is nef and big, hence we can apply Theorem~\ref{Q-fac-abundance}. 
If $C^2<0$, then we can contract such a curve $C_X$ in advance and 
may exclude this case. 
Thus, we can assume that $C^2=0$, i.e. $C_X^2=0$. 
We consider the following equation 
$$K_Y+(2+\alpha)C=f^*(K_X+C_X),$$
where $\alpha$ is the positive integer satisfying $\alpha C=f^*C_X$. 
We apply adjunction to $(X, C_X)$ and $(Y, C)$ (Theorem \ref{adjunction}(2)). 
Then we obtain 
$(K_X+C_X)|_{C_X}\sim_{\Q} 0$ and $(K_Y+C)|_C \sim_{\Q} 0$. 
We deduce that the $\Q$-divisor $C|_C$ is torsion, 
which in turn implies that $C$ is semi-ample 
by Mumford's result (cf. Proposition~\ref{icct}).

If $(X, \Delta)$ is klt, then we can apply a similar argument, 
although we encounter some further technical difficulties in the general case. 
For more details, see Section~\ref{section-klt-abundance}. 

\subsubsection{Log canonical case}

We overview the proof of the log canonical case assuming the klt case. 
After we get the abundance theorem for klt surfaces, 
the main difficulty is to show 
the non-vanishing theorem (Theorem~\ref{non-vanishing}): 
for a log canonical surface $(X, \Delta)$,  
if $K_X+\Delta$ is nef, then $\kappa(X, K_X+\Delta) \geq 0$. 
To this end, we may replace $X$ by its minimal resolution, 
hence assume that $X$ is regular. 
Since we may assume that $\kappa(X, K_X)=-\infty$, 
it follows from the abundance theorem for the regular case that 
a $K_X$-MMP induces a Mori fibre space structure $Z \to B$ 
where $Z$ is the end result. 
Roughly speaking, our idea is to prove that 
the base change of the Mori fibre space to the algebraic closure 
is again a Mori fibre space, 
although the singularities could no longer be log canonical. 
Based on this idea, we extend Fujita's result over algebraically closed fields 
to the case over imperfect fields (cf. Lemma~\ref{semi-positivity}). 
For more details, see Subsection~\ref{Subsection-nonvanishing}.

\begin{rem}
In \cite{BCZ}, 
Birkar--Chen--Zhang independently obtained Theorem~\ref{0abundance} 
under the additional assumption that $(X, \Delta)$ is klt, 
$\Delta$ is a $\Q$-divisor, $S=\Spec\,k$ and $\kappa(X, K_X+\Delta) \geq 0$. 
Their strategy differs from ours.  
\end{rem}

\begin{rem}
In \cite[Corollary 1.2]{Fujino} and \cite[Theorem 1.2]{minimal}, 
Fujino and the author establish the following abundance theorem: 
if $X$ is a projective $\Q$-factorial surface over an algebraically closed field, $\Delta$ is an effective $\R$-divisor 
on $X$ whose coefficients are contained in $[0, 1]$, and $K_X+\Delta$ is nef, 
then $K_X+\Delta$ is semi-ample. 
For some ideas of the proof, we refer to \cite{FT12}. 
The author does not know whether this statement generalises to non-closed fields. 
\end{rem}

\begin{ack}
Part of this work was done whilst 
the author visited National Taiwan University in December 2013 
with the support of the National Center for Theoretical Sciences. 
He would like to thank Professor Jungkai Alfred Chen 
for his generous hospitality. 
The author would like to thank Professors 
Caucher~Birkar, 
Yoshinori~Gongyo, 
Paolo~Cascini, 
J\'anos~Koll\'ar, 
Mircea~Musta\c{t}\u{a}, 
Chenyang~Xu 
for very useful comments and discussions. 
The author also thanks the referee for reading the manuscript carefully and for suggesting several improvements. 
This work was partially supported by JSPS KAKENHI (no. 24224001) and EPSRC. 
\end{ack}

%% file: section2.tex
%%%%%

\section{Preliminaries}

\subsection{Notation}\label{subsection-notation}

Let $k$ be a field. 
We say $X$ is a {\em variety} over $k$ (or $k$-variety) if 
$X$ is an integral scheme which is separated and of finite type over $k$. 
We say $X$ is a {\em curve} over $k$ (resp. a {\em surface} over $k$) 
if $X$ is a variety over $k$ with $\dim X=1$ (resp. $\dim X=2$).
For a scheme $X$, set $X_{\red}$ to be 
the reduced closed subscheme of $X$ 
such that the induced closed immersion $X_{\red} \to X$ is surjective. 

We will not distinguish the notations 
invertible sheaves and divisors. 
For example, 
we will write $L+M$ for 
invertible sheaves $L$ and $M$. 
For a proper scheme $X$ over a field $k$ and a coherent sheaf $F$ on $X$, 
we set 
$$h^i(X, F):=\dim_k H^i(X, F),\,\,\,\chi(X, F):=\sum_{i=0}^{\dim X} (-1)^ih^i(X, F).$$
Note that these numbers depend on $k$. 
For the definition of intersection numbers, 
we refer to \cite[Section 2.2]{MMP}.

We will freely use the notation 
and terminology in \cite{Kollar}. 
In the definition in \cite[Definition~2.8]{Kollar}, 
for a pair $(X, \Delta)$, $\Delta$ is not necessarily effective. 
However, in this paper, we always assume that $\Delta$ is an effective $\R$-divisor.

Let $\Delta$ be an $\R$-divisor on a normal variety over a field. 
For real numbers $a$ and $b$, we write $\Delta \leq a$ (resp. $\Delta \geq b$) 
if, for the prime decomposition $\Delta=\sum_{i\in I}\delta_i \Delta_i$, 
$\delta_i \leq a$ (resp. $\delta_i \geq b$) holds for every $i\in I$. 
For example, $0\leq \Delta \leq 1$ means $0\leq \delta_i \leq 1$ for every $i \in I$. 

For a $\Z$-module $M$, 
we set $M_{\Q}:=M\otimes_{\Z} \Q$ and $M_{\R}:=M\otimes_{\Z} \R$.

%For the definition of dualizing sheaves and canonical divisors, 
%see Subsection~\ref{subsection-dualizing}.

\subsection{Dualising sheaves and purely inseparable base changes}\label{subsection-dualizing}

In this subsection, we summarise basic properties of dualising sheaves and purely inseparable base changes 
(cf. \cite[Subsections~2.2, 2.3]{purely}). 
%For definitions and basic properties of dualising sheaves, 
%see \cite[Definition 1.6]{Kollar} or \cite[Subsection 2.3]{purely}. 

\begin{lem}\label{purely-basic}
Let $k$ be a field of characteristic $p>0$. 
Let $X$ be a normal variety over $k$. 
Then the following assertions hold. 
\begin{enumerate}
\item{The normalisation morphism $Y \to (X\times_k k^{1/p^{\infty}})_{\red}$ is a universal homeomorphism. }
\item{$\rho(X)=\rho(Y)$.}
\item{If $X$ is $\Q$-factorial, then so is $Y$. }
\end{enumerate}
\end{lem}

\begin{proof}
See \cite[Lemma~2.2, Proposition~2.4, Lemma~2.5]{purely}.
\end{proof}

We recall the definition of dualising sheaves and 
collect some basic properties.

\begin{dfn}
Let $k$ be a field. 
Let $X$ be a $d$-dimensional separated scheme of finite type over $k$. 
We set 
$$\omega_{X/k}:=\mathcal{H}^{-d}(f^!k),$$
where $f:X \to \Spec\,k$ is the structure morphism. 
\end{dfn}

A dualising sheaf does not change by enlarging a base field. 

\begin{lem}\label{base-invariance}
Let $k_0 \subset k$ be a field extension with $[k:k_0]<\infty$. 
Let $X$ be a separated scheme of finite type over $k$. 
Note that $X$ is also of finite type over $k_0$. 
Then, there exists an isomorphism
$$\omega_{X/k}\simeq \omega_{X/k_0}.$$ 
\end{lem}

\begin{proof}
See \cite[Lemma~2.7]{purely}. 
\end{proof}

Thanks to Lemma~\ref{base-invariance}, 
we can define canonical divisors $K_X$ independent of a base field.
Our definition is the same as in \cite[Definition 1.6]{Kollar}. 

\begin{dfn}
Let $k$ be a field. 
If $X$ is a normal variety over $k$, 
then it is well-known that $\omega_{X/k}$ is a reflexive sheaf. 
Let $K_X$ be a divisor which satisfies 
$$\MO_X(K_X)\simeq \omega_{X/k}.$$
Such a divisor $K_X$ is called {\em canonical divisor}. 
Note that a canonical divisor is determined up to linear equivalence. 
\end{dfn}

\begin{thm}\label{purely}
Let $k$ be a field of characteristic $p>0$. 
Let $X$ be a normal variety over $k$. 
Let $\nu:Y \to (X\times_k k^{1/p^{\infty}})_{\red}$ be the normalization of $(X\times_k k^{1/p^{\infty}})_{\red}$ and 
set $f:Y \to X$ to be the induced morphism: 
$$f:Y \overset{\nu}\to (X\times_k k^{1/p^{\infty}})_{\red} \to X\times_k k^{1/p^{\infty}} \to X,$$
where the second arrow is taking the reduced structure and 
the third arrow is the first projection. 
If $K_X$ is $\Q$-Cartier, then there exists an effective $\Z$-divisor $E$ on $Y$ such that 
$$K_Y+E=f^*K_X.$$
\end{thm}

\begin{proof}
See \cite[Theorem~4.2]{purely}. 
\end{proof}

\subsection{Some criteria for semi-ampleness}

We collect some known results (at least for algebraically closed fields) on semi-ampleness.

\begin{lem}\label{positive-semiample}
Let $k$ be a field. 
Let $X$ be a projective normal surface over $k$. 
Let $D$ be an effective $\Q$-Cartier $\Q$-divisor and 
let $D=\sum_{i \in I} d_iD_i$ be the irreducible decomposition. 
If $D\cdot D_i>0$ for every $i \in I$, then $D$ is semi-ample. 
\end{lem}

\begin{proof}
We show that $\mathbf B(D)=\emptyset$, 
where $\mathbf B(D)$ is the stable base locus of $D$. 
We obtain $\mathbf B(D) \subset \Supp(D)$. 
On the other hand, by \cite[Theorem~1.3]{Birkar}, 
we obtain $\mathbf B(D) \subset \mathbf B_+(D)=\mathbb E(D).$ 
Thus we see 
$$\mathbf B(D) \subset \Supp(D) \cap \mathbb E(D).$$ 
By the assumption $D\cdot D_i>0$, 
$\Supp(D) \cap \mathbb E(D)$ is empty or zero-dimensional. 
Therefore, \cite[Corollary~1.14]{Fuj83} implies that $D$ is semi-ample. 
\end{proof}

\begin{lem}\label{kappa1}
Let $k$ be a field. 
Let $X$ be a projective normal surface over $k$. 
Let $L$ be a $\Q$-Cartier $\Q$-divisor. 
If $L$ is nef and $\kappa(X, L)=1$, 
then $L$ is semi-ample. 
\end{lem}

\begin{proof}
If $k$ is algebraically closed, 
then the assertion is well-known (cf. \cite[Theorem~4.1]{Fuj84}). 
Thanks to Lemma~\ref{purely-basic}, 
We can reduce the problem to this case by 
taking the normalisation of $(X\times_k \overline{k})_{\red}$. 
\end{proof}

\begin{lem}\label{MFS}
Let $k$ be a field. 
Let $X$ be a projective normal surface over $k$ and 
let $\Delta$ be an effective $\R$-divisor such that $K_X+\Delta$ is $\R$-Cartier. 
Let $L$ be a $\Q$-Cartier $\Q$-divisor on $X$.  
Assume the following conditions. 
\begin{enumerate}
\item{$L$ is nef.}
\item{$L^2=0.$}
\item{There exists a curve $C$ in $X$ with $L\cdot C>0$.}
\item{$L\cdot (K_X+\Delta)<0$.}
\end{enumerate} 
Then, $L$ is semi-ample. 
\end{lem}

\begin{proof}
By replacing $X$ with the minimal resolution, 
we may assume that $X$ is regular and $\Delta=0$. 
By (1) and (3), we see $H^2(X, mL)=0$ for $m\gg 0$. 
Then, by the Riemann--Roch formula 
(cf. \cite[Theorem 2.10]{MMP}), we obtain 
\begin{eqnarray*}
\chi(X, mL)
&=&\frac{1}{2}mL\cdot (mL-K_X)+\chi(X, \MO_X)\\
&=&\frac{1}{2} mL\cdot (-K_X)+\chi(X, \MO_X).\\
\end{eqnarray*} 
Hence, the condition (4) implies $\kappa(X, L)\geq 1$. 
By the condition (2), we obtain $\kappa(X, L)=1$. 
Then, $L$ is semi-ample by Lemma~\ref{kappa1}. 
\end{proof}

\begin{thm}\label{Keel-surface}
Let $k$ be a field of positive characteristic. 
Let $X$ be a projective normal surface over $k$ and 
let $L$ be a nef and big $\Q$-Cartier $\Q$-divisor. 
We define a reduced closed subscheme $\mathbb E(L)$ of $X$ by 
$$\Supp \mathbb E(L)=\bigcup_{L\cdot C=0} C$$
where $C$ runs over curves in $X$. 
If $L|_{\mathbb E(L)}$ is semi-ample, 
then $L$ is semi-ample. 
\end{thm}

\begin{proof}
See \cite[Theorem~0.2]{Keel}.
\end{proof}

%%%%%
\subsection{Contractible criterion}

In the classical minimal model program for surfaces, 
Castelnuovo's criterion plays a crucial role. 
We recall a similar result for singular surfaces. 

\begin{thm}\label{cont}
Let $k$ be a field. 
Let $X$ be a quasi-projective normal surface over $k$. 
Let $C$ be a curve in $X$ proper over $k$. 
Assume that the following conditions hold. 
\begin{itemize}
\item{$X$ is $\Q$-factorial.}
\item{$(K_X+C)\cdot C<0.$}
\item{$C^2<0$.}
\end{itemize}
Then there exists a projective birational morphism $f:X\to Y$ which satisfies the following properties. 
\begin{enumerate}
\item{$Y$ is a quasi-projective normal surface with $f_*\MO_X=\MO_Y$.}
\item{$\Ex(f)=C$.}
\item{Let $L$ be a Cartier divisor on $X$ such that $L\cdot C=0$. 
Then, $mL=f^*L_Y$ for some integer $m>0$ and a Cartier divisor $L_Y$ on $Y$. }
\item{If $X$ is $\Q$-factorial, then $Y$ is $\Q$-factorial.}
\item{If $X$ is projective, then $Y$ is projective and $\rho(Y)=\rho(X)-1$.}
\end{enumerate}
\end{thm}

\begin{proof}
By taking a projective compactification $X \subset \overline X$ such that 
$\overline X$ is regular along $\overline X\setminus X$, 
we can assume that $X$ is projective. 
Then the assertion follows from \cite[Theorem 4.4]{MMP}.
\end{proof}

\subsection{Adjunction}
In this section, we summarise some results on adjunction formula. 

\begin{prop}\label{Fujino-exact}
Let $k$ be a field. 
Let $X$ be a normal $k$-surface and 
let $C$ be a curve in $X$. 
Then, there exists an exact sequence  
$$0\to \mathcal T \to (\omega_X\otimes_{\MO_X} \MO_X(C))|_C \to \omega_C \to 0$$
where $\mathcal T$ is a skyscraper sheaf.
\end{prop}

\begin{proof}
See the proof of \cite[Lemma~4.4]{Fujino}.
\end{proof}

\begin{lem}\label{adjunction-lemma}
Let $k$ be a field. 
Let $X$ be a normal $k$-surface and 
let $C$ be a proper $k$-curve in $X$. 
Let $L$ be a Cartier divisor on $X$. 
If $H^1(C, \MO_X(L)|_C)\neq 0$, then 
$$H^0(C, (\MO_X(r(K_X+C-L))|_C)\neq 0$$ 
for every $r\in\mathbb Z_{>0}$. 
%where we define $F^{*}:=\mathcal Hom_{\MO_X}(F, \MO_X)$ for a coherent sheaf $F$ on $X$. 
\end{lem}

\begin{proof}
Set $\omega_X(C-L):=\omega_X \otimes_{\MO_X} \MO_X(C-L)$.
By Proposition~\ref{Fujino-exact}, 
we obtain an exact sequence: 
$$0\to \mathcal T\otimes_{\MO_C} \MO_X(-L)|_C \to (\omega_X(C-L))|_C \to 
\omega_C\otimes_{\MO_C} (\MO_X(-L)|_C) \to 0.$$
Since $\mathcal T$ is a skyscraper sheaf, we have 
$H^1(C, \mathcal T\otimes_{\MO_C} \MO_X(-L)|_C)=0.$ 
Since 
$$H^0(C, \omega_C\otimes_{\MO_C} (\MO_X(-L)|_C)) \simeq H^1(C, \MO_X(L)|_C)^*\neq 0,$$ 
there is a nonzero element 
$$\xi \in H^0(C, \omega_X(C-L)|_C)$$ 
whose image $\xi' \in H^0(C, \omega_C\otimes_{\MO_C} (\MO_X(-L)|_C)$ is non-zero as well. 
Thus there exists a map 
$$\mathcal O_{C} \to \omega_X(C-L)|_C$$
such that this is injective on some non-empty open set. 
Tensoring $(\omega_X(C-L))|_C$ one by one, 
we obtain a sequence of maps 
$$\mathcal O_{C} 
\to (\omega_X(C-L))|_C
\to (\omega_X(C-L))^{\otimes 2}|_C
\to \cdots
\to (\omega_X(C-L))^{\otimes r}|_C
$$
which are injective on some non-empty open set. 
On the other hand, 
there is a natural map 
$$(\omega_X(C-L))^{\otimes r}|_C \to ((\omega_X(C-L))^{\otimes r})^{**}|_C=\MO_X(r(K_X+C-L))|_C,$$
which is bijective on some non-empty open set. 
Combining these maps, we have a map 
$$\mathcal O_{C} \to \MO_X(r(K_X+C-L))|_C,$$
which is injective on some non-empty open set. 
Thus, the kernel $K$ of this map is a torsion subsheaf of $\mathcal O_C$. 
Then, we have $K=0$. 
Therefore, we obtain an injection $\mathcal O_{C} \hookrightarrow \MO_X(r(K_X+C-L))|_C.$ 
This implies $H^0(C, \MO_X(r(K_X+C-L))|_C) \neq 0.$
\end{proof}

Using Lemma~\ref{adjunction-lemma}, 
we obtain the following result.

\begin{thm}\label{adjunction}
Let $k$ be a field. 
Let $X$ be a normal $k$-surface and 
let $C$ be a proper $k$-curve in $X$ such that $r(K_X+C)$ is Cartier for some positive integer $r$. 
\begin{enumerate}
\item{Assume $(K_X+C)\cdot C<0$. 
If $L$ is a Cartier divisor on $X$ such that $L\cdot C=0$, 
then $L|_C\simeq \MO_C$. }
\item{If $(K_X+C)\cdot C=0$, then 
there exists $s\in r\mathbb Z_{>0}$ such that 
$\mathcal{O}_X(s(K_X+C))|_C \simeq \mathcal O_C$.}
\end{enumerate}
\end{thm}

\begin{proof}
(1) 
It is enough to show $H^0(C, L|_C)\neq 0$. 
Since $(K_X+C)\cdot C<0$ implies $H^0(C, \MO_X(r(K_X+C-L))|_C)=0$, 
we obtain $H^1(C, L|_C)=0$ by Lemma~\ref{adjunction-lemma}. 
In particular, we see $H^1(C, \MO_C)=0$. 
By the Riemann--Roch formula, we obtain 
$$h^0(C, L|_C)=\chi(C, L|_C)=\chi(C, \MO_C)=h^0(C, \MO_C)\neq 0.$$

(2) 
By replacing $r$ with $2r$, we can assume $r\geq 2$. 
Assume $H^1(C, \MO_C)\neq 0$. 
Then, we can apply Lemma~\ref{adjunction-lemma} for $L:=0$ and 
we obtain 
$$H^0(C, r(K_X+C)|_C)\neq 0.$$
Thus, we may assume $H^1(C, \MO_C)=0$. 

Assume $H^1(C, r(K_X+C)|_C)\neq 0$. 
Then, we can apply Lemma~\ref{adjunction-lemma} for $L:=r(K_X+C),$ 
and we obtain 
$$H^0(C, -(r-1)r(K_X+C)|_C)\neq 0.$$
Thus, we may assume $H^1(C, r(K_X+C)|_C)=0.$ 

It is enough to show that $H^0(C, r(K_X+C)|_C)\neq 0$, assuming
$$H^1(C, \MO_C)=H^1(C, r(K_X+C)|_C)=0.$$
This follows from the Riemann--Roch formula: 
$$h^0(C, r(K_X+C)|_C)=\chi(C, r(K_X+C)|_C)=\chi(C, \MO_C)=h^0(C, \MO_C)\neq 0.$$
\end{proof}

\subsection{Semiample perturbation}

In this subsection, 
we show Lemma \ref{perturb-leq1}. 
Let us start with a consequence of Noether's normalisation theorem.

\begin{lem}\label{perturb-leq0}
Let $k$ be an infinite field. 
Let $f:X \to W$ be a projective surjective $k$-morphism of $k$-varieties. 
Then there exists a non-empty open subset $W^0 \subset W$ such that 
if $X^0:=f^{-1}(W^0)$, then the induced morphism 
$f^0:X^0 \to W^0$ factors as follows:  
\[
f^0:X^0 \xrightarrow{\pi} W^0 \times_k \mathbb P^r_{k} \xrightarrow{{\rm pr_1}}  W^0, 
\]
where $r=\dim X-\dim W$, 
$\pi$ is a finite surjective morphism and ${\rm pr_1}$ denotes the first projection. 
\end{lem}

\begin{proof}
Since $K(W)$ is an infinite field, Noether's normalisation theorem deduces that 
there is a finite surjective morphism: 
\[
\widetilde{\pi}:X \times_W \Spec\,K(W) \to \mathbb P^r_{K(W)},
\]
where $r:=\dim X \times_W \Spec\,K(W)= \dim X-\dim W$. 
By killing denominators, we obtain what we want. 
\end{proof}

\begin{lem}\label{perturb-leq1}
Let $k$ be an infinite field. 
Let $X$ be a projective normal $k$-variety and 
let $D$ be a semi-ample $\Q$-Cartier $\Q$-divisor. 
If $\epsilon\in \Q_{>0}$ and $x_1, \cdots, x_{\ell} \in X$, then 
there exists an effective $\Q$-Cartier $\Q$-divisor $D'$ 
such that $D' \sim_{\Q} D$, $0\leq D' \leq \epsilon$, and that $x_i \not\in\Supp D'$ for every $i$. 
\end{lem}

\begin{proof}
By replacing $D$ with a large multiple, 
we may assume that $\epsilon=1$ and that 
$D$ is a base point free Cartier divisor whose linear system induces a morphism 
$$f:X \to W$$
with $f_*\MO_X=\MO_W$ and $D=f^*\MO_W(1)$, where $\MO_W(1)$ is very ample. 
Since $f$ is projective, we can apply 
Lemma \ref{perturb-leq0}, hence the following hold. 
\begin{itemize}
\item{$W^0 \subset W$ is a non-empty open subset and set $X^0:=f^{-1}(W^0)$.}
\item{The induced morphism $f|_{X^0}:X^0 \to W^0$ can be written by 
$$f|_{X^0}:X^0 \overset{\pi} \to W^0 \times_k \mathbb P^r_{k} \overset{{\rm pr_1}}\to W^0,$$
where $\pi$ is a finite surjective morphism. }
\end{itemize}
Since $D$ is base point free, 
it suffices to show that there exists $\nu \in \Z_{>0}$ such that 
every general hyperplane section $H \sim \MO_W(1)$ satisfies 
$$0\leq f^*H \leq \nu.$$
Let 
$$\pi:X^0 \overset{\pi_i}\to V \overset{\pi_s}\to W^0 \times_k \mathbb P^r_{k}$$
be the factorisation corresponding to the separable closure of 
$K(W^0 \times_k \mathbb P^r_{k})$ in $K(X^0)$. 
Take a positive integer $e$ such that $K(X^0) \subset K(V)^{1/p^e}$. 
Since the absolute Frobenius morphism $F^e_V:V \to V$ is the same as the normalisation of $V$ 
in $K(V)^{1/p^e}$, 
the absolute Frobenius $F^e_V$ factors through $X^0$: 
$$F^e_V:V \to X^0 \overset{\pi_i}\to V.$$
Then, we can check that for every prime divisor $P$ on $W^0$, 
we obtain 
$$\pi^*(P) \leq p^e\deg(\pi_s)=:\nu.$$
Here, note that $\pi^*(P)$ is defined 
as the closure of the pull-back of the restriction of $P$ to the regular locus of $W^0$. 
For a general hyperplane section $H \subset W$, 
its pull-back $f^*H$ is equal to the closure of $(f|_{X^0})^*(H|_{W^0})$. 
By \cite[Corollary 1]{Seidenberg}, 
a general hyperplane section $H \sim \MO_W(1)$ is reduced. 
Thus, also ${\rm pr}_1^*(H|_{W^0})$ is reduced, that is, ${\rm pr}_1^*(H|_{W^0}) \leq 1$. 
Therefore, we obtain 
$$\pi^*{\rm pr}_1^*(H|_{W^0}) \leq p^e\deg(\pi_s)=\nu.$$
This implies $0\leq f^*H \leq \nu.$
We are done. 
\end{proof}

%% file: section3.tex
%%%%%%%%%%%
\section{Abundance theorem for klt surfaces}\label{section-klt-abundance}

In this section, we prove a special case of the abundance theorem (Theorem~\ref{sep-closed}). 
This theorem implies the abundance theorem for klt surfaces (Theorem~\ref{klt-abundance}). 
First we give a criterion for semi-ampleness.

\begin{prop}\label{icct}
Let $k$ be a field of characteristic $p>0$. 
Let $X$ be a projective normal surface over $k$. 
Let $E$ be an effective $\Z$-divisor on $X$ and 
let $E=\sum_{j \in J} e_j E_j$ be the prime decomposition. 
Assume the following conditions. 
\begin{enumerate}
\item{$K_X$ and $E$ are $\Q$-Cartier.}
\item{$K_X\cdot E_j=E\cdot E_j=0$ for every $j \in J$.}
\item{There exists $s\in \Z_{>0}$ such that 
$sE$ is Cartier and that $\MO_X(sE)|_E \simeq \MO_E$. }
\end{enumerate}
Then, $E$ is semi-ample. 
\end{prop}

\begin{proof}
The proof consists of four steps.

\setcounter{step}{0}

\begin{step}\label{s1-icct}
In order to prove the assertion of Proposition \ref{icct}, 
we may assume that 
\begin{enumerate}
\item[(a)] $k$ is separably closed, and
\item[(b)] $X$ is regular. 
\end{enumerate} 
\end{step}

\begin{proof}[Proof of Step \ref{s1-icct}] 
After replacing $k$ by $H^0(X, \MO_X)$, 
we may assume that $X$ is geometrically irreducible (Remark \ref{r-icct}). 
Taking the base change of $X$ to the separable closure of $k$, 
we may assume (a). 

Let $h:X' \to X$ be a resolution of singularities. 
Let $n$ be a positive integer such that $nE$ is Cartier. 
Replacing $X$ and $E$ with $X'$ and $h^*(nE)$ respectively, 
we may assume (b). 
This completes the proof of Step \ref{s1-icct}. 
\end{proof}

\begin{step}\label{s2-icct}
The assertion of Proposition \ref{icct} holds 
if $k$ is algebraically closed. 
\end{step}

\begin{proof}[Proof of Step \ref{s2-icct}] 
By Step \ref{s1-icct}, we may assume that $X$ is smooth over $k$. 
Then the assertion of Proposition \ref{icct} 
follows directly from \cite[Lemma in page 682]{Masek}. 
This completes the proof of Step \ref{s2-icct}. 
\end{proof}

\begin{step}\label{step-irreducible}
In order to prove the assertion of Proposition \ref{icct}, 
we may assume that 
\begin{enumerate}
\item[(a)] $k$ is separably closed,
\item[(c)] $X$ is $\Q$-factorial, and 
\item[(d)] $E$ is irreducible.  
\end{enumerate} 
\end{step}

\begin{proof}[Proof of Step \ref{step-irreducible}] 
By Step~\ref{s1-icct}, 
we may assume that (a) and (b) hold. 
In particular, also (c) holds. 
Clearly we may assume that $E$ is connected. 

Let us  reduce the proof of Proposition \ref{icct} to the case when $E$ is irreducible. 
Assume that $E$ is not irreducible. 
Fix an irreducible component $E_1$ of $E$. 
Since $E$ is connected and $E\cdot E_1=0$, we obtain $E_1^2<0$. 
This implies $(K_X+E_1)\cdot E_1<0$. 
By Theorem~\ref{cont}, there exists a birational morphism 
$${\rm cont}_{E_1}=\pi:X \to Z$$
to a projective $\Q$-factorial surface $Z$ such that $\Ex(\pi)=E_1$. 

For the time being, we assume that $Z$ and $E':=\pi_*E$ satisfies the same properties as (1), (2), and (3). 
If $E'$ is irreducible, then we are done. 
Otherwise, we can apply the same argument as above. 
Since the number of the irreducible components of $E'$ 
is strictly less than the one of $E$, this repeating procedure will terminate. 

Therefore, it suffices to show that $Z$ and $E'$ satisfy (1), (2), and (3). 
Since $Z$ is $\Q$-factorial, the property (1) holds automatically. 
By $K_X\cdot E_1=0$, we can write 
$$K_X=\pi^*K_Z.$$ 
By $E\cdot E_1=0$, %and Theorem~\ref{cont}(3), 
we obtain $E=\pi^*\pi_*E=\pi^*(E')$. 
Set $E_j':=\pi_*E_j$ for $j \in J \setminus \{1\}$. 
Thus $K_Z$ and $E'$ satisfy the property (2) by 
$$K_Z\cdot E_j'=K_Z\cdot \pi_*E_j=\pi^*K_Z\cdot E_j=K_X \cdot E_j=0$$
and 
$$E'\cdot E_j'=E'\cdot \pi_*E_j=\pi^*E'\cdot E_j=E \cdot E_j=0.$$
We show that $E'$ satisfies the condition (3), that is, $\MO_Z(tE')|_{E'} \simeq \MO_{E'}$ for some integer $t>0$. 
Set 
$$E\overset{\alpha}\to E'' \overset{\beta}\to E'$$ 
to be the Stein factorisation. 
Since every fibre of $E \to E'$ is connected and $k$ is separably closed, 
$E'' \to E'$ is a universally homeomorphism. 
Thus it suffices to show that 
$$\beta^*(\MO_Z(tE')|_{E'}) \simeq \MO_{E''}$$
for some $t>0$. 
Fix $t>0$ such that $tE$ and $tE'$ are Cartier and that 
$\MO_X(tE)|_{E} \simeq \MO_E$. 
Then, we obtain the required isomorphism $\beta^*(\MO_Z(tE')|_{E'}) \simeq \MO_{E''}$ by 
$$\alpha_*(\MO_X(tE)|_{E}) \simeq \alpha_*\MO_E \simeq \MO_{E''}$$
and 
$$\alpha_*(\MO_X(tE)|_{E}) \simeq \alpha_*(\pi^*\MO_Z(tE')|_{E}) 
\simeq \alpha_*(\alpha^*\beta^*(\MO_Z(tE')|_{E'})) \simeq \beta^*(\MO_Z(tE')|_{E'}).$$
This completes the proof of Step~\ref{step-irreducible}. 
\end{proof}

\begin{step}\label{step-to-AlgClosed}
The assertion of Proposition \ref{icct} holds 
without any additional assumptions. 
\end{step}

\begin{proof}[Proof of Step \ref{step-to-AlgClosed}]
Thanks to Step~\ref{step-irreducible}, 
we may assume that (a), (c) and (d) hold. 
Let $Y$ be the normalisation of $(X\times_k \overline{k})_{\red}$ and 
let 
$$f:Y \to X$$ 
be the induced morphism. 
By Theorem~\ref{purely}, we obtain 
$$K_Y+D=f^*K_X$$
for some effective $\Z$-divisor $D$ on $Y$. 
Note that $Y$ is $\Q$-factorial by Lemma~\ref{purely-basic}(3). 
Set $E_Y:=f^*E$. 

We now show that $Y$ and $E_Y$ satisfy the properties (1), (2), and (3). 
Since $Y$ is $\Q$-factorial, (1) holds automatically. 
Since 
$$\MO_E \simeq \MO_X(sE)|_{E}$$
for some $s \in \Z_{>0}$, taking the pull-backs to $E_Y$, 
we obtain 
$$\MO_{E_Y} \simeq (f^*\MO_X(sE))|_{E_Y} \simeq \MO_Y(sE_Y)|_{E_Y}.$$ 
Thus $Y$ and $E_Y$ satisfy (3). 
We see $E_Y^2=0$ by $E^2=0$. 
In particular $E_Y$ is nef. 
Thus, we obtain $D\cdot E_Y\geq 0.$ 
If $K_Y \cdot E_Y<0$, then $E_Y$ is semi-ample by Lemma~\ref{MFS}. 
Therefore, we may assume that $K_Y\cdot E_Y\geq 0$. 
Since $(K_Y+D)\cdot E_Y=f^*K_X \cdot E_Y=0,$ 
the inequalities $K_Y\cdot E_Y\geq 0$ and $D\cdot E_Y\geq 0$ imply 
$$K_Y\cdot E_Y=D\cdot E_Y=0.$$
Thus (2) holds. 
Thus, $Y$ and $E_Y$ satisfies the conditions (1), (2), and (3). 

Therefore, Step \ref{s2-icct} implies that $E_Y$ is semi-ample. 
Hence, also $E$ is semi-ample. 
This completes the proof of Step \ref{step-to-AlgClosed}. 
\end{proof}

Step~\ref{step-to-AlgClosed} completes the proof of Proposition \ref{icct}. 
\end{proof}

\begin{rem}\label{r-icct}
Let $k$ be a field and let $X$ be a projective normal variety over $k$. 
Then $k':=H^0(X, \MO_X)$ is a field. Moreover, $k'$ gives the Stein factorisation 
of the structure morphism: 
\[
X \xrightarrow{\alpha'} \Spec\,k' \to \Spec\,k.
\]
In particular, $X$ is geometrically connected over $k'$. 
Moreover, it follows from \cite[Lemma 2.2(1)]{purely} that $X$ is geometrically irreducible over $k'$. 
\end{rem}

We prove the main result of this section.

\begin{thm}\label{sep-closed}
Let $k$ be a separably closed field of characteristic $p>0$. 
Let $X$ be a projective normal $\Q$-factorial surface over $k$ and 
let $\Delta$ be a $\Q$-divisor with $0\leq \Delta \leq 1$. 
Assume that the following condition $(*)$ holds. 
\begin{enumerate}
\item[$(*)$]{For a curve $C$ in $X$, if $(K_X+\Delta)\cdot C=0$, then $C^2 \geq 0$. }
\end{enumerate}
If $K_X+\Delta$ is nef, then $K_X+\Delta$ is semi-ample. 
\end{thm}

\begin{proof}
Since $k$ is separably closed, $X$ is geometrically irreducible. 
Let $Y$ be the normalisation of $(X\times_k \overline k)_{\red}$ and 
set
$$f:Y \to X$$
to be the induced morphism. 
Then, $Y$ is $\Q$-factorial (Lemma~\ref{purely-basic}(3)). 
Thanks to Theorem~\ref{purely}, we obtain 
$$K_{Y}+E+f^*\Delta=f^*(K_X+\Delta)$$
for some effective $\Z$-divisor $E$ on $Y$. 
Let 
$$E+f^*\Delta=\sum_i C_i$$
be the decomposition into the connected components and 
let 
$$C_i=\sum_j c_{ij}C_{ij}$$
be the irreducible decomposition.

\setcounter{step}{0}

\begin{step}\label{step-non-negative}
In order to prove Theorem \ref{sep-closed}, 
we may assume that 
\begin{enumerate}
\item[(a)] $(K_X+\Delta)^2=0$, and 
\item[(b)] for every irreducible component $C$ of $\Supp E \cup \Supp f^*\Delta$, we obtain 
$$f^*(K_X+\Delta)\cdot C=0.$$
In  particular, $C^2\geq 0$. 
\end{enumerate}
\end{step}

\begin{proof}[Proof of Step \ref{step-non-negative}]
For a rational number $0<\epsilon < 1$, we consider the following equation: 
$$K_Y+E+f^*\Delta=\epsilon K_Y+(\epsilon(E+f^*\Delta)+(1-\epsilon)f^*(K_X+\Delta)).$$
Note that, for $0<\epsilon \ll 1$, the latter divisor 
$$\epsilon(E+f^*\Delta)+(1-\epsilon)f^*(K_X+\Delta)$$ 
is nef by $(*)$. 
We fix such a small rational number $0<\epsilon \ll 1$. 
If this divisor is big, then we obtain 
$$K_Y+E+f^*\Delta=\epsilon K_Y+(\epsilon (E+f^*\Delta)+(1-\epsilon)f^*(K_X+\Delta))\sim_{\mathbb Q}\epsilon(K_Y+\Delta_Y)$$
where $\Delta_Y$ is a $\Q$-divisor with $0 \leq \Delta_Y \leq 1$. 
%whose coefficients are contained in $[0, 1]$. 
Then, we can apply \cite[Theorem~1.2]{minimal}. 
Thus, we may assume $(K_X+\Delta)^2=0$. Hence (a) holds. 
Moreover, we may assume that, for every irreducible component $C$ of $\Supp E \cup \Supp f^*\Delta$, 
$$f^*(K_X+\Delta)\cdot C=0.$$ 
Hence (b) holds. 
Then the condition $(*)$ implies $C^2\geq 0$. 
This completes the proof of Step \ref{step-non-negative}.  
\end{proof}

\begin{step}\label{step-perturb}
Assume that (a) and (b) of Step \ref{step-non-negative} holds. 
Fix an arbitrary index $i$. 
Then the following hold. 
\begin{enumerate}
\item{If $C_i$ is reducible, then $C_i$ is semi-ample. }
\item{Assume that $C_i$ is irreducible. 
Let $C_i=c_{i0}C_{i0}$ where $c_{i0}\in\mathbb Q_{>0}$ and 
$C_{i0}$ is the prime divisor. 
\begin{enumerate}
\item[(2a)]{If $C_i^2>0$, then $C_{i}$ is semi-ample.}
\item[(2b)]{If $C_i^2=0$, then $C_i$ is semi-ample or $0\leq c_{i0}\leq 1$. }
\end{enumerate}}
\end{enumerate}
\end{step}

\begin{proof}[Proof of Step \ref{step-perturb}]
(1) Assume $C_i$ is reducible. 
Then $C_{ij}^2 \geq 0$ implies $C_i \cdot C_{ij} >0$ for every $j$. 
We see that $C_i$ is semi-ample by Lemma~\ref{positive-semiample}. 

(2a) 
The assertion follows from Lemma~\ref{positive-semiample}. 

(2b) 
Assuming that $c_{i0}>1$, it suffices to show that $C_i$ is semi-ample. 
Set $C_{i0, X}:=f(C_{i0})$ to be the curve. 
Since $(K_X+\Delta)\cdot C_{i0, X}=0$ (cf. Step~\ref{step-non-negative}) and $C_{i0, X}^2=0$, 
we obtain $K_X \cdot C_{i0, X}\leq 0$. 
If $K_X\cdot C_{i0, X}<0$, then $C_{i0, X}$ is semi-ample by Lemma~\ref{MFS}. 
Thus, we may assume 
$$K_X\cdot C_{i0, X}=\Delta\cdot C_{i0, X}=0.$$ 
This implies $(K_Y+E+f^*\Delta)\cdot C_{i0}=0$. 
Since $C_{i0}^2=0$, we obtain $K_Y \cdot C_{i0}\leq 0$. 
If $K_Y\cdot C_{i0}<0$, then $C_{i0}$ is semi-ample by Lemma~\ref{MFS}. 
Thus we can assume 
$$K_Y\cdot C_{i0}=(E+f^*\Delta)\cdot C_{i0}=0.$$
This implies that $C_{i0}$ (resp. $C_{i0, X}$) 
is a connected component of $\Supp(E+f^*\Delta)$ (resp. $\Supp \Delta$). 
Thus, for sufficiently divisible $s\in \Z_{>0}$, we obtain 
\begin{eqnarray*}
\MO_X(s(K_X+\Delta))|_{C_{i0, X}} &\simeq& \MO_X(s(K_X+\delta C_{i0, X}))|_{C_{i0, X}},\,\,\,\,{\rm and}\\
\MO_Y(s(K_Y+E+f^*\Delta))|_{C_{i0}} &\simeq& \MO_Y(s(K_Y+c_{i0}C_{i0}))|_{C_{i0}},
\end{eqnarray*}
where $0\leq \delta \leq 1$ is the coefficient of $C_{i0, X}$ in $\Delta$. 
Since $(K_X+C_{i0, X})\cdot C_{i0, X}=0$ and $(K_Y+C_{i0})\cdot C_{i0}=0$, 
we can apply Theorem \ref{adjunction}(2) and obtain 
\begin{eqnarray*}
\MO_X(s(K_X+C_{i0, X}))|_{C_{i0, X}} &\simeq& \MO_{C_{i0, X}},\,\,\,\,{\rm and}\\
\MO_Y(s(K_Y+C_{i0}))|_{C_{i0}} &\simeq& \MO_{C_{i0}}
\end{eqnarray*}
for sufficiently divisible $s\in \Z_{>0}$. 
Summarising above, we have 
\begin{eqnarray*}
\MO_X(s(K_X+\Delta))|_{C_{i0, X}} &\simeq & \MO_X(s(-(1-\delta)C_{i0, X}))|_{C_{i0, X}}
,\\
\MO_Y(s(K_Y+E+f^*\Delta))|_{C_{i0}} &\simeq & \MO_Y(s((c_{i0}-1)C_{i0}))|_{C_{i0}}. 
\end{eqnarray*}
Since $\MO_Y(s(K_Y+E+f^*\Delta)) \simeq f^*\MO_X(s(K_X+\Delta))$, 
we obtain 
\begin{eqnarray*}
&&\MO_Y(s((c_{i0}-1)C_{i0}+(1-\delta)f^*C_{i0, X})))|_{C_{i0}}\\
&\simeq& \MO_Y(s(K_Y+E+f^*\Delta-f^*(K_X+\Delta)))|_{C_{i0}}\\
&\simeq& \MO_{C_{i0}}.
\end{eqnarray*}

Thus, for some sufficiently divisible $s' \in \Z_{>0}$, we obtain 
$$\MO_Y(s'C_{i0})|_{C_{i0}} \simeq \MO_{C_{i0}}.$$
Since $K_Y\cdot C_{i0}=C_{i0}^2=0$, 
we can apply Lemma~\ref{icct} and $C_{i0}$ is semi-ample. 
This completes the proof of Step \ref{step-perturb}
\end{proof}

\begin{step}\label{s3-sep-closed}
The assertion of Theorem \ref{sep-closed} holds without any additional assumptions. 
\end{step}

\begin{proof}[Proof of Step \ref{s3-sep-closed}]
By Step \ref{step-non-negative}, 
the problem is reduced to the case when 
(a) and (b) of Step \ref{step-non-negative} hold. 
Then (1) and (2) of Step~\ref{step-perturb} hold.  
Therefore, Lemma~\ref{perturb-leq1} enables us to find 
a $\Q$-divisor $\Delta_Y$ on $Y$ such that $0 \leq \Delta_Y \leq 1$ 
and 
$$E+f^*\Delta=\sum_i C_i\sim_{\Q} \Delta_Y.$$ 
Then the divisor $K_Y+E+f^*\Delta\sim_{\Q}K_Y+\Delta_Y$ is semi-ample by \cite[Theorem~1.2]{minimal}. 
This completes the proof of Step \ref{s3-sep-closed}. 
\end{proof}
Step \ref{s3-sep-closed} completes the proof of Theorem \ref{sep-closed}. 
\end{proof}

As consequences of Theorem~\ref{sep-closed}, 
we obtain the abundance theorem for the klt case. 

\begin{thm}\label{regular-abundance}
Let $k$ be a separably closed field of characteristic $p>0$. 
Let $X$ be a projective normal $\Q$-factorial surface over $k$ and 
let $\Delta$ be a $\Q$-divisor with $\llcorner \Delta \lrcorner=0$. 
If $K_X+\Delta$ is nef, then $K_X+\Delta$ is semi-ample. 
\end{thm}

\begin{proof}
We reduce the proof to the case when 
the condition $(*)$ in Theorem~\ref{sep-closed} holds. 
Assume that there is a curve $C$ such that $(K_X+\Delta)\cdot C=0$ and that $C^2<0$. 
Since $0\leq \Delta <1$, we obtain 
$$(K_X+C) \cdot C< (K_X+\Delta)\cdot C=0.$$
Then, by Theorem~\ref{cont}, 
we obtain a contraction $f:X\to X'$ of $C$ to a $\Q$-factorial surface. 
We repeat this procedure and this will terminate. 
Thus we may assume that $(*)$ in Theorem~\ref{sep-closed} holds, 
hence the assertion follows from Theorem~\ref{sep-closed}. 
\end{proof}

\begin{thm}\label{klt-abundance}
Let $k$ be a field of characteristic $p>0$. 
Let $(X, \Delta)$ be a projective klt surface over $k$, 
where $\Delta$ is an effective $\Q$-divisor. 
If $K_X+\Delta$ is nef, then $K_X+\Delta$ is semi-ample. 
\end{thm}

\begin{proof}
Taking the base change to the separable closure after replacing $k$ by $H^0(X, \MO_X)$ (Remark \ref{r-icct}), 
we may assume that $k$ is separably closed. 
Taking the minimal resolution of $X$, we may assume that $X$ is regular. 
Then the assertion follows from Theorem~\ref{regular-abundance}. 
\end{proof}

%% file: section4.tex
\section{Abundance theorem for log canonical surfaces}\label{Section-lc-abundance}

In this section, we show the abundance theorem for log canonical surfaces (Theorem~\ref{abundance}), 
that is, 
for a projective log canonical surface $(X, \Delta)$ with a $\Q$-divisor $\Delta$, 
if $K_X+\Delta$ is nef, then it is semi-ample. 
Subsection~\ref{Subsection-nonvanishing} is devoted to show that $\kappa(X, K_X+\Delta)\geq 0$. 
In Subsection~\ref{Subsection-abundance}, we prove that $K_X+\Delta$ is semi-ample 
for each case $\kappa(X, K_X+\Delta)=0, 1,$ and $2$. 
In Subsection~\ref{Subsection-relative-abundance}, we generalise our result to the relative settings.

\subsection{Non-vanishing theorem}\label{Subsection-nonvanishing}

The goal of this subsection is to show Theorem~\ref{non-vanishing}. 
A rough overview of this subsection is as follows. 
To show the non-vanishing theorem (Theorem~\ref{non-vanishing}), 
we can assume that $\kappa(X, K_X)=-\infty$. 
Thanks to the abundance theorem for the klt case (Theorem~\ref{klt-abundance}), 
we may assume that $(X\times_k \overline k)_{\red}$ is a ruled surface. 
The rational case (resp. the irrational case) is treated in Proposition~\ref{nonvanishing-rational} 
(resp. Proposition~\ref{nonvanishing-irrational}).

\begin{prop}\label{nonvanishing-rational}
Let $k$ be a separably closed field of characteristic $p>0$. 
Let $(X, \Delta)$ be a projective log canonical surface over $k$, 
where $\Delta$ is a $\Q$-divisor. 
Assume the following conditions. 
\begin{itemize}
\item{$X$ is regular.}
\item{$(X\times_k \overline k)_{\red}$ is a rational surface.}
\end{itemize}
If $K_X+\Delta$ is nef, then $\kappa(X, K_X+\Delta) \geq 0$. 
\end{prop}

\begin{proof}
%We can replace the minimal resolution and we may assume that $X$ is regular. 
Let $Y$ be the minimal resolution of the normalisation $(X\times_k \overline k)_{\red}^N$ 
of $(X\times_k \overline k)_{\red}$: 
$$f:Y \to (X\times_k \overline k)_{\red}^N \to (X\times_k \overline k)_{\red} 
\to X\times_k \overline k \to X.$$ 
By our assumption, $Y$ is a smooth rational surface. 

If $K_X+\Delta\equiv 0$, then the pull-back to $Y$ is torsion, hence also $K_X+\Delta$ is torsion by Lemma~\ref{purely-basic}. 
Thus we may assume that $K_X+\Delta \not\equiv 0$.  
Then, by Theorem~\ref{purely}, we obtain 
$$K_Y+\Delta_Y=f^*(K_X+\Delta)$$
for some effective $\Q$-divisor $\Delta_Y$. 
Fix $m_0\in \Z_{>0}$ such that $m_0(K_Y+\Delta_Y)$ is Cartier. 
Then, we can check that $H^2(Y, mm_0(K_Y+\Delta_Y))=0$ for every large integer $m\gg 0$. 
Therefore, by the Riemann--Roch formula (cf. \cite[Theorem 2.10]{MMP}), we obtain 
\begin{eqnarray*}
&&h^0(Y, mm_0(K_Y+\Delta_Y))\\
&\geq& \chi(Y, \MO_Y)+\frac 1 2 mm_0(K_Y+\Delta_Y)\cdot(mm_0(K_Y+\Delta_Y)-K_Y)\\
&=&1+\frac 1 2 mm_0(K_Y+\Delta_Y)\cdot(\Delta_Y+(mm_0-1)(K_Y+\Delta_Y))\\
&\geq &1,
\end{eqnarray*}
where the above equality follows from the fact that $Y$ is rational. 
\end{proof}

To show Lemma~\ref{semi-positivity}, 
we establish auxiliary results on Mori fibre spaces: 
Lemma~\ref{pic-1-fano} and Lemma~\ref{genus0-MFS}. 

\begin{lem}\label{pic-1-fano}
Let $k$ be a separably closed field of characteristic $p>0$. 
Let $X$ be a projective klt $k$-surface such that $-K_X$ is ample and that $\rho(X)=1$. 
Then $(X\times_k \overline k)_{\red}$ is a rational surface. 
\end{lem}

\begin{proof}
If $k=\overline{\mathbb F}_p$, then the assertion is well-known. 
Thus we can assume that $k\neq \overline{\mathbb F}_p$. 
Since $\overline{\mathbb F}_p \subsetneq k$, 
we get $\overline k\neq \overline{\mathbb F}_p$. 

Assume that $(X\times_k \overline k)_{\red}$ is not a rational surface and 
let us derive a contradiction. 
Set $Y$ to be the normalisation of $(X\times_k \overline k)_{\red}$. 
Then, by Lemma~\ref{purely-basic}, 
$Y$ is a projective normal $\Q$-factorial surface such that $\rho(Y)=1$. 
Moreover,  by Theorem~\ref{purely} and $\rho(Y)=1$, $-K_Y$ is ample. 
Let $f:Z \to Y$ be the minimal resolution. 
We obtain $K_Z+E=f^*K_Y$ for some effective $\Q$-divisor $E$ on $Z$. 
Since $Z$ is not rational and $K_Z$ is not pseudo-effective, 
by running $K_Z$-MMP, we see that $Z$ has a ruled surface structure $\pi:Z \to B$. 
Since $\overline k\neq \overline{\mathbb F}_p$ and $Z$ is not rational, we can apply \cite[Theorem~3.20]{minimal}. 
Therefore, every $f$-exceptional curves are $\pi$-vertical. 
Thus, $Z \to B$ factors through $Y$. 
In particular, the induced morphism $Y \to B$ is a surjection. 
However, this contradicts $\rho(Y)=1$. 
\end{proof}

\begin{lem}\label{genus0-MFS}
Let $k$ be a separably closed field of characteristic $p>0$. 
Let $g:Z \to B$ be a surjective morphism 
from a regular projective $k$-surface to a regular projective $k$-curve. 
Assume that $g$ is a $K_Z$-Mori fibre space structure, that is, 
$g_*\MO_Z=\MO_B$, $\rho(Z/B)=1$, and $-K_Z$ is $g$-ample. 
Consider the following commutative diagram: 
$$\begin{CD}
Z @<<< Z\times_k \overline k @<<< (Z\times_k \overline k)_{\red} @<<< (Z\times_k \overline k)_{\red}^N=:W\\
@VVg V @VVg\times_k {\overline k} V @VV(g\times_k {\overline k})_{\red} V 
@VV(g \times_k {\overline k})_{\red}^N V\\
B @<<< B\times_k \overline k @<<< (B\times_k \overline k)_{\red} @<<< (B\times_k \overline k)_{\red}^N=:\Gamma',
\end{CD}$$
where $V^N$ denotes the normalisation of a variety $V$. 
Let 
$$(g \times_k {\overline k})_{\red}^N:W \overset{h}\to \Gamma \overset{\gamma}\to \Gamma'$$
be the Stein factorisation of $(g \times_k {\overline k})_{\red}^N$.
Then the following assertions hold. 
\begin{enumerate}
\item{$W$ is $\Q$-factorial.}
\item{$\gamma:\Gamma \to \Gamma'$ is a universal homeomorphism.}
\item{A general fibre of $h:W \to \Gamma$ is $\mathbb P^1_{\overline k}$.}
\item{If $(Z\times_k \overline k)_{\red}$ is not rational, then 
$H^1(B, \MO_B)\neq 0$ and neither $\Gamma$ nor $\Gamma'$ is a rational curve. }
\end{enumerate}
\end{lem}

\begin{proof}
(1) 
The assertion follows from Lemma~\ref{purely-basic}.

(2) 
Since $g_*\MO_Z=\MO_B$ and the base change $(-)\times_k \overline{k}$ is flat, 
we obtain $(g\times_k {\overline k})_*\MO_{Z\times_k \overline{k}}=\MO_{B\times_k \overline{k}}$. 
Moreover, the horizontal arrows 
$$(Z\times_k \overline k)_{\red}^N \to Z\times_k \overline k,\,\,\, (B\times_k \overline k)_{\red}^N \to B\times_k \overline k$$ 
are universal homeomorphisms (Lemma~\ref{purely-basic}). 
Thus every fibre of $(g \times_k {\overline k})_{\red}^N$ is connected. 
Therefore $\gamma:\Gamma \to \Gamma'=(B\times_k \overline k)_{\red}^N$ 
is a finite morphism whose fibres are connected. 
In particular, $\gamma$ is a universal homeomorphism. 

(3) 
Let $f:W \to Z$ be the induced morphism. 
By Theorem~\ref{purely}, 
we can find an effective $\Z$-divisor $D$ on $W$ such that 
$$K_W+D=f^*K_Z.$$
Let $F_g$ (resp. $F_h$) be a general fibre of $g$ (resp. $h$). 
By (2), $K_Z\cdot F_g<0$ implies $K_W\cdot F_h<0$. 
Thus, we see $(K_W+F_h)\cdot F_h<0$. It follows from \cite[Theorem 3.19(1)]{minimal} 
that $F_h \simeq \mathbb P^1_{\overline{k}}$.

(4) 
First, we show that 
neither $\Gamma$ nor $\Gamma'$ is a rational curve. 
By (3), there exists a non-empty open subset 
$\Gamma^0 \subset \Gamma$ such that 
$h^{-1}(\Gamma^0) \simeq \Gamma^0 \times_{\overline k} \mathbb P^1_{\overline k}$. 
Since $(Z\times_k \overline k)_{\red}$ is not rational, 
$\Gamma$ is not a rational curve. 
By (2), $\gamma:\Gamma \to \Gamma'$ is a finite surjective purely inseparable morphism. 
Therefore $\Gamma'$ is not a rational curve.

Second we show $H^1(B, \MO_B)\neq 0$. 
Assume that  $H^1(B, \MO_B)=0$ and let us derive a contradiction. 
We obtain $\deg(K_B)<0$. 
%By \cite[Lemma~10.6]{Kollar}, $B$ is isomorphic to $\mathbb P^1_K$ or a regular conic in $\mathbb P^2_K$, 
%where $K:=H^0(B, \MO_B)$. 
By Theorem~\ref{purely}, 
$\Gamma':=(B\times_k \overline{k})_{\red}^N$ satisfies $\deg(K_{\Gamma'})<0$. 
Thus $\Gamma' \simeq \mathbb P^1_{\overline{k}}$. 
However, we have already shown $\Gamma' \not\simeq \mathbb P^1_{\overline{k}}$, 
which is a contradiction. 
\end{proof}

\begin{rem}\label{h0-vector-space}
Lemma~\ref{genus0-MFS}(4) states that $H^1(B, \MO_B) \neq 0$ 
for a projective regular curve $B$. 
This implies 
$$\chi(B, \MO_B) =\dim_k H^0(B, \MO_B)-\dim_k H^1(B, \MO_B)\leq 0$$ 
as follows. 

Let $B \to \Spec\,k_B \to \Spec\,k$ be the Stein factorisation of 
the structure morphism. 
Then, we obtain $H^0(B, \MO_B) \simeq k_B$. 
On the other hand, $H^1(B, \MO_B)$ has a $k_B$-vector space structure and we obtain 
\[
\dim_{k}H^i(B, \MO_B)=[k_B: k]\dim_{k_B}H^i(B, \MO_B).
\]
Therefore, $H^1(B, \MO_B) \neq 0$ implies 
\begin{eqnarray*}
&&\dim_k H^0(B, \MO_B)-\dim_k H^1(B, \MO_B)\\
&=&[k_B: k](\dim_{k_B} H^0(B, \MO_B)-\dim_{k_B} H^1(B, \MO_B))\\
&=&[k_B: k](1-\dim_{k_B} H^1(B, \MO_B)) \leq 0.
\end{eqnarray*}
\end{rem}

The following lemma is a key result to show the non-vanishing theorem for the irrational case 
(Proposition~\ref{nonvanishing-irrational}). 

\begin{lem}\label{semi-positivity}
Let $k$ be a separably closed field of characteristic $p>0$. 
Let $g:Z \to B$ be a surjective morphism 
from a regular projective surface $Z$ to a regular projective curve $B$. 
Let $\Delta_Z$ be a $\Q$-divisor on $Z$ with $0\leq \Delta_Z \leq 1$. 
Assume the following conditions. 
\begin{itemize}
\item{$g$ is a $K_Z$-Mori fibre space structure, that is, 
$g_*\MO_Z=\MO_B$, $\rho(Z/B)=1$, and $-K_Z$ is $g$-ample. }
\item{$K_Z+\Delta_Z$ is $g$-nef. }
\item{$(Z\times_k \overline k)_{\red}$ is not rational.}
\end{itemize}
Then $\kappa(Z, K_Z+\Delta_Z) \geq 0$. 
\end{lem}

\begin{proof}

\setcounter{step}{0}

\begin{step}\label{s1-semi-positivity}
In order to prove Lemma \ref{semi-positivity}, 
we may assume that 
\begin{enumerate}
\item $\Supp\,\Delta_Z$ contains no $g$-vertical prime divisor, and 
\item there exists a $\Q$-Cartier $\Q$-divisor $L_B$ on $B$ such that 
$\deg\,L_B \leq 0$ and 
$K_Z+\Delta_Z \sim_{\Q} g^*L_B$. 
\end{enumerate}
\end{step}

\begin{proof}[Proof of Step \ref{s1-semi-positivity}]
Take a general fibre $F_g$ of $g$. 
Note that $F_g$ is irreducible. 
By dropping the $g$-vertical components of $\Delta_Z$, we may assume (1). 
Moreover, by reducing the coefficients of $\Delta_Z$, 
the problem is reduced to the case when $(K_Z+\Delta_Z) \cdot F_g=0.$ 
It follows from $\rho(Z/B)=1$ that $K_Z+\Delta_Z \equiv_g 0.$ 
By \cite[Theorem 4.4]{MMP}, we obtain $K_Z+\Delta_Z=g^*L_B$ for some 
$\Q$-Cartier $\Q$-divisor $L_B$ on $B$. 
If $L_B$ is ample, then there is nothing to show. 
Thus, we may assume that $\deg(L_B) \leq 0$. 
Hence, (2) holds. 
This completes the proof of Step \ref{s1-semi-positivity}. 
\end{proof}

From now on, we assume that (1) and (2) of Step \ref{s1-semi-positivity} hold. 

\begin{step}\label{s2-semi-positivity}
The assertion of Lemma \ref{semi-positivity} holds 
if there exists a $g$-horizontal curve $C$ such that $C^2 \leq 0$. 
\end{step}

\begin{proof}[Proof of Step \ref{s2-semi-positivity}]
We have 
$$0\geq g^*L_B\cdot C=(K_Z+\Delta_Z)\cdot C$$ 
$$\geq (K_Z+C)\cdot C=2(h^1(C, \MO_C)-h^0(C, \MO_C))\geq 0,$$
where the second equality follows by  \cite[Corollary 2.8]{MMP} and the last inequality 
holds by Lemma~\ref{genus0-MFS}(4) and Remark~\ref{h0-vector-space}. 
Thus all the above inequalities are equalities. 
In particular, we obtain $\deg\,L_B=0$. 

We now treat the case when $C^2 <0$. 
For a sufficiently divisible integer $s\in \Z_{>0}$, 
we get  
$$g^*(sL_B)|_C \simeq \MO_Z(s(K_Z+\Delta_Z))|_C \simeq \MO_Z(s(K_Z+C))|_C,$$
where the first isomorphism follows from Step \ref{s1-semi-positivity}(2) 
and the second isomorphism holds by $C^2<0$ and the equation 
$(K_Z+\Delta_Z)\cdot C = (K_Z+C)\cdot C$. 
By Theorem \ref{adjunction}(2), this is torsion. 
Therefore, $L_B$ is also torsion, hence its pull-back $g^*L_B \simeq K_Z+\Delta_Z$ is semi-ample. 
This completes the proof for the case when $C^2<0$. 

Hence, we may assume that $C^2 = 0$ and $D^2 \geq 0$ for any curve $D$ on $Z$. 
Since $K_Z+\Delta_Z=g^*L_B \equiv 0$, 
Theorem~\ref{sep-closed} implies that $K_Z+\Delta_Z$ is semi-ample. 
In any case, we see that $K_Z+\Delta_Z$ is semi-ample. 
This completes the proof of Step \ref{s2-semi-positivity}. 
\end{proof}

\begin{step}\label{s3-semi-positivity}
The assertion of Lemma \ref{semi-positivity} holds 
without any additional assumptions. 
\end{step}

\begin{proof}[Proof of Step \ref{s3-semi-positivity}]
By Step \ref{s2-semi-positivity}, 
we can assume that every $g$-horizontal curve has positive self-intersection number. 
It follows from Lemma \ref{positive-semiample} that every curve on $Z$ is semi-ample. 

Consider the following commutative diagram: 
$$\begin{CD}
Z @<<< Z\times_k \overline k @<<< (Z\times_k \overline k)_{\red} @<<< (Z\times_k \overline k)_{\red}^N=:W\\
@VVg V @VVg\times_k {\overline k} V @VV(g\times_k {\overline k})_{\red} V 
@VV(g \times_k {\overline k})_{\red}^N V\\
B @<<< B\times_k \overline k @<<< (B\times_k \overline k)_{\red} @<<< (B\times_k \overline k)_{\red}^N=:\Gamma',
\end{CD}$$
where $Y^N$ denotes the normalisation of a variety $Y$. 
Let 
$$(g \times_k {\overline k})_{\red}^N:W \overset{h}\to \Gamma \overset{\gamma}\to \Gamma'$$
be the Stein factorisation of $(g \times_k {\overline k})_{\red}^N$ and 
we obtain the following commutative diagram 
$$\begin{CD}
Z @<f << W\\
@VVg V @VVh V\\
B @<\rho << \Gamma. 
\end{CD}$$
By Theorem~\ref{purely}, we can find an effective $\Q$-divisor $\Delta_W$ on $W$ such that 
$$K_W+\Delta_W=f^*(K_Z+\Delta)=f^*g^*L_B=h^*\rho^*L_B.$$
We see that $\Delta_W$ is semi-ample. 
%We obtain $K_W \cdot F_h<0$ for a general fiber $F_h$ of $h$. 
Let $\mu:V \to W$ 
be the minimal resolution and set 
$$q:V \overset{\mu}\to W \overset{h}\to \Gamma.$$ 
We can find an effective $\mu$-exceptional $\Q$-divisor $E_V$ on $V$ such that 
$$K_V+E_V=\mu^*K_W,$$
$$K_V+E_V+\mu^*\Delta_W=\mu^*(K_W+\Delta_W).$$
Clearly, every $\mu$-exceptional divisor is $q$-vertical. 
Therefore, the divisor 
$$K_V+\mu^*\Delta_W=\mu^*(K_W+\Delta_W)-E_V=\mu^*h^*\rho^*L_B-E_V$$
satisfies $(K_V+\mu^*\Delta_W)\cdot F_q=0$ for any general fibre $F_q$ of $q:V \to \Gamma$. 

Since $\Delta_W$ is semi-ample, 
Lemma \ref{perturb-leq1} enables us to find a $\Q$-divisor $\Delta_V$ on $V$ 
such that $\Delta_V \sim_{\Q} \mu^*\Delta_W$ and $0 \leq \Delta_V \leq 1$. 
Therefore, $\kappa(V, K_V+\mu^*\Delta_W) \geq 0$ by \cite[Theorem 2.2]{Fuj84}. 
Note that we can apply \cite[Theorem 2.2]{Fuj84} 
because $\Gamma$ is irrational by Lemma~\ref{genus0-MFS}(4). 
In particular, we obtain 
$$\kappa(Z, K_Z+\Delta_Z)=\kappa(W, K_W+\Delta_W)=\kappa(V, K_V+\mu^*\Delta_W+E_V)\geq 0.$$
This completes the proof of Step \ref{s3-semi-positivity}. 
\end{proof}
Step \ref{s3-semi-positivity} completes the proof of Lemma \ref{semi-positivity}. 
\end{proof}

We prove the non-vanishing theorem for the irrational case. 

\begin{prop}\label{nonvanishing-irrational}
Let $k$ be a separably closed field of characteristic $p>0$. 
Let $(X, \Delta)$ be a projective log canonical surface over $k$, 
where $\Delta$ is a $\Q$-divisor. 
Assume the following conditions. 
\begin{itemize}
\item{$X$ is regular.}
\item{$(X\times_k \overline k)_{\red}$ is not a rational surface.}
\item{$\kappa(X, K_X)=-\infty$. }
\end{itemize}
If $K_X+\Delta$ is nef, then $\kappa(X, K_X+\Delta) \geq 0$. 
\end{prop}

\begin{proof}
We run a $K_X$-MMP with scaling $\Delta$. 
In other words, there exists a sequence 
of birational morphisms of regular projective surfaces: 
\[
X=:X_0 \xrightarrow{f_0} X_1 \xrightarrow{f_1} \cdots 
\xrightarrow{f_{\ell-1}} X_{\ell}
\]
and a $K_{X_{\ell}}$-Mori fibre space 
\[
f_{\ell}:X_{\ell} \to X_{\ell+1}
\]
such that 
\begin{enumerate}
\item $-K_{X_i}$ is $f_i$-ample for any $i \in \{0, \cdots, \ell\}$, 
\item $\rho(X_i/X_{i+1})=1$ for any $i \in \{0, \cdots, \ell\}$, 
and  
\item $K_{X_i}+\lambda_i \Delta_i \equiv_{f_i} 0$ for any $i \in \{0, \cdots, \ell\}$, where 
$\Delta_i$ is the push-forward of $\Delta$ on $X_i$ and 
$\lambda_i$ is defined by 
\[
\lambda_i:=\inf\{\mu \in \R_{\geq 0}\,|\, K_{X_i}+\mu \Delta_i \text{ is nef}\}.
\] 
\end{enumerate}
It is well-known that 
\begin{itemize}
\item $1 \geq \lambda_0 \geq \lambda_1 \geq \cdots \geq \lambda_{\ell} \geq 0$, and 
\item 
for any $i$, $\lambda_i$ is a rational number. 
\end{itemize}
In particular, we have that 
\begin{itemize}
\item $\kappa(X, K_X+\Delta) \geq \kappa(X_{\ell}, K_{X_{\ell}}+\lambda_{\ell}\Delta_{\ell})$, and 
\item $K_{X_{\ell}}+\lambda_{\ell}\Delta_{\ell}$ is nef over $X_{\ell+1}$. 
\end{itemize}
By Lemma~\ref{pic-1-fano}, we see that $\dim X_{\ell+1}=1$. 
Then we obtain $\kappa(X_{\ell}, K_{X_{\ell}}+\lambda_{\ell}\Delta_{\ell}) \geq 0$ by Lemma~\ref{semi-positivity}. 
Therefore, we get $\kappa(X, K_X+\Delta) \geq 0$, as desired. 
\end{proof}

We show the main result of this subsection.

\begin{thm}\label{non-vanishing}
Let $k$ be a field of characteristic $p>0$. 
Let $(X, \Delta)$ be a projective log canonical surface over $k$, 
where $\Delta$ is a $\Q$-divisor. 
If $K_X+\Delta$ is pseudo-effective, then $\kappa(X, K_X+\Delta) \geq 0$. 
\end{thm}

\begin{proof}
By \cite[Theorem 1.1]{MMP}, 
we may assume that $K_X+\Delta$ is nef. 
Taking the base change to the separable closure of $k$ 
after replacing $k$ by $H^0(X, \MO_X)$ (Remark \ref{r-icct}), 
we may assume that $k$ is separably closed. 
By replacing $X$ with its minimal resolution, 
we may assume that $X$ is regular. 

If $\kappa(X, K_X)\geq 0$, then there is nothing to show. 
If $\kappa(X, K_X)=-\infty$, 
then the assertion follows from 
Proposition~\ref{nonvanishing-rational} and Proposition~\ref{nonvanishing-irrational}. 
\end{proof}

\subsection{Abundance theorem}\label{Subsection-abundance}

In this subsection, 
we show the abundance theorem with $\Q$-coefficients (Theorem~\ref{abundance}). 
In Proposition~\ref{kappa-0} (resp. Proposition~\ref{kappa-2}), 
we treat the case $\kappa(X, K_X+\Delta)=0$ (resp. $\kappa(X, K_X+\Delta)=2$).

\begin{prop}\label{kappa-0}
Let $k$ be a separably closed field of characteristic $p>0$. 
Let $X$ be a projective normal $\Q$-factorial surface over $k$ and 
let $\Delta$ be a $\Q$-divisor on $X$ with $0\leq \Delta \leq 1$. 
If $K_X+\Delta$ is nef and $\kappa(X, K_X+\Delta)=0$, then $K_X+\Delta$ is semi-ample. 
\end{prop}

\begin{proof} 
Since $\kappa(X, K_X+\Delta)=0$, 
we obtain 
$$K_X+\Delta \sim_{\mathbb Q} D$$
for some effective $\Q$-divisor $D$. 
We assume $D \neq 0$ and let us derive a contradiction. 
Let 
$$D=\sum_{i\in I} d_iD_i$$
be the irreducible decomposition with $d_i \in \Q_{>0}$. 

\setcounter{step}{0}
\begin{step}\label{step-negative}
In order to prove Proposition \ref{kappa-0}, 
we may assume that 
if $C$ is a curve on $X$ such that $(K_X+\Delta)\cdot C=0$, 
then 
\begin{enumerate}
\item{$C^2 \geq 0$, or }
\item{$C \subset \Supp(\llcorner \Delta \lrcorner)$ and 
$C$ is a connected component of $\Supp \Delta$. }
\end{enumerate} 
\end{step}

\begin{proof}[Proof of Step \ref{step-negative}]
Let $C$ be a curve in $X$ such that 
$(K_X+\Delta)\cdot C=0$ and that $C$ satisfies none of (1) and (2). 
We obtain $C^2<0$. 
Hence, we get 
$$(K_X+C)\cdot C \leq (K_X+\Delta)\cdot C=0.$$
If this inequality is an equality, then (2) holds. 
Thus we obtain $(K_X+C)\cdot C<0$. 
By Theorem~\ref{cont}, 
we obtain a birational morphism $h:X \to Y$ to a projective $\Q$-factorial surface $Y$ 
such that $\Ex(h)=C$. 
We can check that $Y$ and $\Delta_Y:=h_*\Delta$ satisfy the same properties as $X$ and $\Delta$, 
i.e.  $0\leq \Delta_Y\leq 1$ and $K_Y+\Delta_Y$ is nef. 
Moreover, if $K_Y+\Delta_Y$ is semi-ample, 
then also its pullback $K_X+\Delta=h^*(K_X+\Delta)$ is semi-ample. 
Thus the problem is reduced to the one of $(Y, \Delta_Y)$. 
If there exists a curve $C'$ on $Y$, with $(K_Y+\Delta_Y)\cdot C'=0$, 
which satisfies none (1) and (2), 
then we can repeat the same procedure as above. 
This procedure will terminate because the Picard number strictly drops:  $\rho(Y)=\rho(X)-1$. 
Therefore, we may assume that every curve $C$ on $X$, with $(K_X+\Delta)\cdot C=0$, 
satisfies (1) or (2). 
This completes the proof of Step \ref{step-negative}. 
\end{proof}

From now on, we assume that 
if $C$ is a curve on $X$ such that $(K_X+\Delta)\cdot C=0$, 
then (1) or (2) of Step \ref{step-negative} holds.

\begin{step}\label{s2-kappa-0}
Any connected component of $D$ is irreducible. 
\end{step}

\begin{proof}[Proof of Step \ref{s2-kappa-0}]
Assume that $D_{i_1} \cap D_{i_2} \neq\emptyset$ 
for some irreducible components $D_{i_1} \neq D_{i_2}$ of $\Supp\,D$. 
Since $D$ is nef and 
$$D^2=(K_X+\Delta)\cdot D=0,$$ 
we obtain $D\cdot D_{i_a}=(K_X+\Delta)\cdot D_{i_a}=0$ for any $a\in \{1, 2\}$. 
Since $D_{i_1} \cdot D_{i_2}>0$, we obtain $D_{i_a}^2<0$ for any $a\in \{1, 2\}$. 
By Step~\ref{step-negative}, each $D_{i_a}$ is a connected component of $\Supp \Delta$. 
This contradicts $D_{i_1} \cap D_{i_2} \neq\emptyset$. 
This completes the proof of Step \ref{s2-kappa-0}. 
\end{proof}

\begin{step}\label{s3-kappa-0}
It holds that 
$$K_X\cdot D_i=\Delta\cdot D_i=D\cdot D_i=D_i^2=0$$
for every $i\in I$. 
\end{step}

\begin{proof}[Proof of Step \ref{s3-kappa-0}]
Since $K_X+\Delta \sim_{\Q} D$ is nef, we obtain 
$$(K_X+\Delta)\cdot D_i=D\cdot D_i=D_i^2=0$$
for every $i\in I$, where the second equation follows from Step \ref{s2-kappa-0}. 
If $K_X\cdot D_i < 0$ for some $i \in I$, then we obtain 
$\kappa(X, D) \geq \kappa (X, D_i)\geq 1$ 
by Lemma~\ref{MFS}. 
This contradicts 
$\kappa(X, K_X+\Delta)=\kappa(X, D)=0.$ 
Thus we get $K_X\cdot D_i \geq 0$ for every $i \in I$.
Since $D_i$ is nef, we obtain $\Delta\cdot D_i \geq 0$. 
Thus, $(K_X+\Delta)\cdot D_i=0$ implies 
$$K_X\cdot D_i=\Delta\cdot D_i=0.$$
This completes the proof of Step \ref{s3-kappa-0}. 
\end{proof}

\begin{step}\label{s4-kappa-0}
The assertion of Proposition \ref{kappa-0} holds 
without any additional assumptions. 
\end{step}

\begin{proof}[Proof of Step \ref{s4-kappa-0}]
Fix $i \in I$. 
It suffices to show that $D_i$ is semi-ample. 
By Step \ref{s2-kappa-0}, we can write 
$$\Delta=\delta D_i+\Delta'$$ 
for some rational number $\delta$ and effective $\Q$-divisor $\Delta'$ such that $0\leq \delta \leq 1$ and $D_i \cap \Supp\,\Delta'=\emptyset$. 
For any sufficiently divisible $s \in \Z_{>0}$, 
we have that 
\begin{eqnarray*}
0
&\sim& \MO_X(s(K_X+D_i))|_{D_i}\\
&=&\MO_X(s(K_X+D_i+(\Delta-\delta D_i)))|_{D_i}\\
&=&\MO_X(s(K_X+\Delta+(1-\delta)D_i))|_{D_i}\\
&\sim& \MO_X(s(D+(1-\delta)D_i))|_{D_i}\\
&=&\MO_X(s(d_i+(1-\delta))D_i)|_{D_i},
\end{eqnarray*}
where the first linear equivalence follows from Step \ref{s3-kappa-0} 
and Theorem \ref{adjunction}(2). 
By $d_i+(1-\delta)\geq d_i>0$, 
we obtain $\MO_X(rD_i)|_{D_i} \simeq \MO_{D_i}$ for a sufficiently divisible integer $r$. 
Therefore, $D_i$ is semi-ample 
by the equation
$$K_X\cdot D_i=D_i^2=0$$
and Proposition~\ref{icct}. 
This completes the proof of Step \ref{s4-kappa-0}. 
\end{proof}
Step \ref{s4-kappa-0} completes the proof of Proposition \ref{kappa-0}. 
\end{proof}

Although the following argument is the same as in \cite[Proposition~3.29]{minimal}, 
we give a proof for the sake of the completeness. 

\begin{prop}\label{kappa-2}
Let $k$ be a field of characteristic $p>0$. 
Let $X$ be a projective normal $\Q$-factorial surface and 
let $\Delta$ be a $\Q$-divisor with $0\leq \Delta \leq 1$. 
If $K_X+\Delta$ is nef and big, then $K_X+\Delta$ is semi-ample. 
\end{prop}

\begin{proof}
Set 
$$E:=\bigcup\limits_{C\cdot (K_X+\Delta)=0}C=C_1 \cup \dots \cup C_r.$$ 
For any curve $C$ with $C \subset E$, we have 
\[
(K_X+C)\cdot C \leq (K_X+\Delta)\cdot C=0.
\]

\setcounter{step}{0}
\begin{step}\label{s1-kappa-2}
In order to prove Proposition \ref{kappa-2}, 
we may assume that if $C$ is a curve in $X$ such that $C \subset E$, then $(K_X+C)\cdot C=0$. 
\end{step}

\begin{proof}[Proof of Step \ref{s1-kappa-2}]
Let $C$ be a curve in $X$ such that 
$C \subset E$ and $(K_X+C)\cdot C<0$. 
Then we obtain $(K_X+\Delta) \cdot C=0$. 
Since $K_X+\Delta$ is big, we can write $K_X+\Delta=A+D$, 
where $A$ is an ample $\Q$-divisor and $D$ is an effective $\Q$-divisor. 
Therefore, it holds that $C^2<0$. 
By Theorem~\ref{cont}, 
there exists a birational morphism 
$f:X \to Y$ to a projective $\Q$-factorial surface $Y$ such that 
$\Ex(f)=C$. 
Set $\Delta_Y:=f_*(\Delta)$. 
Since $K_X+\Delta=f^*(K_Y+\Delta_Y)$ and $Y$ is $\Q$-factorial, 
if we can prove that $K_Y+\Delta_Y$ is semi-ample, 
then also $K_X+\Delta$ is semi-ample. 
We can repeat this procedure and we obtain the required reduction. 
This completes the proof of Step \ref{s1-kappa-2}. 
\end{proof}

From now on, we assume that 
if $C$ is a curve in $X$ such that $C \subset E$, then $(K_X+C)\cdot C=0$.

\begin{step}\label{s2-kappa-2}
If $C$ is a curve in $X$ such that $C \subset E$, 
then $C$ is a connected component of $E$ and 
$\MO_X(s(K_X+\Delta))|_C=\MO_X(s(K_X+C))|_C$ 
for any sufficiently divisible $s\in\Z_{>0}$. 
\end{step}

\begin{proof}[Proof of Step \ref{s2-kappa-2}] 
Take a curve $C$ in $X$ with $C \subset E$. 
In particular, we obtain $(K_X+\Delta) \cdot C=0$. 
By Step \ref{s1-kappa-2}, we have $(K_X+C)\cdot C=0$. 
Therefore, it holds that  
\[
\Delta \cdot C = -K_X \cdot C = C^2. 
\]
We can write $\Delta=aC +\Delta'$ for some $0 \leq a \leq 1$ and an effective $\Q$-divisor $\Delta'$ such that $C \not\subset \Supp\,\Delta'$. We obtain the following inequalities: 
\[
\Delta \cdot C = (aC + \Delta') \cdot C \geq aC^2 \geq C^2. 
\]
Thanks to the equation $\Delta \cdot C = C^2$, 
both the inequality must be equalities. 
This implies that $C \subset \Supp(\llcorner\Delta\lrcorner)$ and 
$C$ is disjoint from any other irreducible component of $\Delta$. 
This completes the proof of Step \ref{s2-kappa-2}. 
\end{proof}

It follows from Step \ref{s1-kappa-2} and Theorem \ref{adjunction}(2) 
that $\MO_X(s(K_X+\Delta))|_C$ is base point free for any curve $C$ with $C \subset E$ and any sufficiently divisible $s \in \Z_{>0}$. 
By Step 2, $\MO_X(s(K_X+\Delta))|_E$ is base point free. 
Thanks to Theorem \ref{Keel-surface}, it holds that $K_X+\Delta$ is semi-ample, as desired.
\end{proof}

We prove the main theorem of this subsection. 

\begin{thm}[Abundance theorem]\label{abundance}
Let $k$ be a field of characteristic $p>0$. 
Let $(X, \Delta)$ be a projective log canonical surface over $k$, where $\Delta$ is a $\Q$-divisor. 
If $K_X+\Delta$ is nef, then $K_X+\Delta$ is semi-ample. 
\end{thm}

\begin{proof}
Taking the base change to the separable closure after replacing $k$ by $H^0(X, \MO_X)$ (Remark \ref{r-icct}), 
we may assume that $k$ is separably closed. 
Set $\kappa:=\kappa(X, K_X+\Delta)$. 
By Theorem~\ref{non-vanishing}, we obtain $\kappa\geq 0$. 
If $\kappa=0$ (resp. $\kappa=1$, resp. $\kappa=2$), 
then $K_X+\Delta$ is semi-ample by 
Proposition~\ref{kappa-0} (resp. Lemma~\ref{kappa1}, resp. Proposition~\ref{kappa-2}). 
\end{proof}

\subsection{Relative version}\label{Subsection-relative-abundance}

The purpose of this section is to prove Theorem~\ref{relative-abundance}, 
which generalises our abundance theorem (Theorem~\ref{abundance}) to the relative case. 
We start with a lemma to compare nef divisors to relative nef divisors. 

\begin{lem}\label{nef-bound}
Let $k$ be a field of characteristic $p>0$. 
Let $X$ be a projective normal surface over $k$ and 
let $\Delta$ be an effective $\R$-divisor such that $K_X+\Delta$ is $\R$-Cartier. 
Let $f:X \to S$ be a morphism to a projective $k$-scheme $S$. 
Let $A_S$ be an ample invertible sheaf on $S$. 
If $K_X+\Delta$ is $f$-nef, then $K_X+\Delta+mf^*A_S$ is nef for some $m\in \Z_{>0}$. 
\end{lem}

\begin{proof}
Taking the base change to the separable closure after replacing $k$ by $H^0(X, \MO_X)$ (Remark \ref{r-icct}), 
we can assume that $k$ is separably closed. 
Taking the Stein factorisation of $f:X\to S$, 
we may assume $f_*\MO_X=\MO_S$. 
Moreover, by taking the Stein factorisation of the structure morphism $S \to \Spec\,k$, 
we can assume that $H^0(X, \MO_X)=H^0(S, \MO_S)=k.$ 
%=\alpha_*\MO_X=\MO_{\Spec\,k}$, 
%where $\alpha:X \to \Spec\,k$ is the structure morphism. 

Let $Y$ be the normalisation of $(X\times_k \overline{k})_{\red}$ and 
let $g:Y \to X$ be the induced morphism. 
By Theorem~\ref{purely}, 
we obtain 
$$K_Y+\Delta_Y=g^*(K_X+\Delta)$$
for some effective $\R$-divisor $\Delta_Y$ on $Y$. 
Note that $g_S^*A_S$ is ample, where $g_S:S\times_k \overline k \to S$. 
If $K_Y+\Delta_Y+mf'^*g_S^*A_S$ is nef for the induced morphism 
$f':Y \to S \times_k \overline k$, then so is $K_X+\Delta+mf^*(A_S)$. 
Thus, we may assume that $k$ is algebraically closed. 
Then, the assertion follows from \cite[Theorem~3.13]{minimal}. 
\end{proof}

We prove the main result of this section. 

\begin{thm}\label{relative-abundance}
Let $k$ be a field of characteristic $p>0$. 
Let $(X, \Delta)$ be a log canonical surface over $k$, where $\Delta$ is a $\Q$-divisor.  
Let $f:X \to S$ be a projective morphism to a separated scheme $S$ of finite type over $k$. 
If $K_X+\Delta$ is $f$-nef, then $K_X+\Delta$ is $f$-semi-ample. 
\end{thm}

\begin{proof}
In order to use \cite[Theorem~1]{Perturbation}, 
we reduce the proof to the case when $k$ is an $F$-finite field containing $\overline{\mathbb F}_p$. 
For this, first we take the base change to 
the composite field $k\overline{\mathbb F}_p$, 
where $k\overline{\mathbb F}_p$ is the minimum subfield of $\overline k$ 
that contains $k$ and $\overline{\mathbb F}_p$. 
By replacing $k$ with $k\overline{\mathbb F}_p$, 
we may assume that $\overline{\mathbb F}_p \subset k$. 
Second, take an intermediate field $\overline{\mathbb F}_p \subset k_1 \subset k$ 
such that $k_1$ is finitely generated over $\overline{\mathbb F}_p$ and 
that all schemes, morphisms and divisors are defined over $k_1$, 
i.e. there exists $X_1$ is a scheme of finite type over $k_1$ with $X_1 \times_{k_1} k \simeq X$ etc. 
By replacing $k$ with $k_1$, 
we can assume that $k$ is an $F$-finite field containing $\overline{\mathbb F}_p$. 

\medskip

We reduce the proof to the case when $X$ and $S$ are projective. 
Since the problem is local on $S$, we may assume that $S$ is affine. 
We can find projective compactifications $S \subset \overline S$ and 
$X \subset \overline X$, that is, 
there exists a commutative diagram 
$$\begin{CD}
X @>>> \overline X\\
@VVf V @VV\overline fV\\
S @>>> \overline S
\end{CD}$$
such that $\overline X$ and $\overline S$ are projective and 
each horizontal arrow is an open immersion. 
By replacing $\overline X$ with a resolution along $\overline X \setminus X$, 
we may assume that $\overline X$ is regular along $\overline X \setminus X$. 
Taking more blowups, we may assume that the support of 
the closure $\Supp(\overline{\Delta})$ is regular 
at every point contained in $\overline X \setminus X$. 
In particular, $(\overline X, \overline{\Delta})$ is log canonical. 
Note that $(\overline X, \overline{\Delta})$ may not be $\overline f$-nef. 
However, by running a $(K_{\overline{X}}+\overline{\Delta})$-MMP over $\overline S$ (\cite[Theorem 1.1]{MMP}), 
the end result $(\overline{X}', \overline{\Delta}')$ is log canonical and 
a minimal model over $\overline S$. 
Therefore, by replacing $(X, \Delta) \to S$ with $(\overline{X}', \overline{\Delta}') \to \overline S$, 
we can assume that $X$ and $S$ are projective.

Fix an ample invertible sheaf $A_S$ on $S$. 
By Lemma~\ref{nef-bound}, 
$K_X+\Delta+mf^*A_S$ is nef for some $m\in \Z_{>0}$. 
We can find $\Delta'\sim_{\Q} \Delta+mf^*A_S$ such that $(X, \Delta')$ is log canonical 
by \cite[Theorem~1]{Perturbation}. 
Thus, by Theorem~\ref{abundance}, 
$K_X+\Delta'$ is semi-ample. 
In particular, $K_X+\Delta$ is $f$-semi-ample. 
\end{proof}

\subsection{Abundance theorem with $\mathbb{R}$-coefficients}

In this subsection, we generalise our relative log canonical abundance theorem 
(Theorem~\ref{relative-abundance}) to the case of $\R$-coefficients (Theorem~\ref{R-abundance}). 
The main strategy is to use Shokurov polytope, 
which is the same as the case when $k$ is algebraically closed. 

However, there are some differences as follows. 
Let $X$ be a projective normal (geometrically connected) variety over a field $k$ and 
let $Y$ be the normalization of $(X\times_k \overline{k})_{\red}$. 
Set $f:Y \to X$ to be the induced homomorphism. 
In the proof of our main result (Theorem~\ref{R-abundance}), 
we consider the following $\Z$-bilinear homomorphism: 
$$M\times N \to \Z,\,\,\, (D, C) \mapsto f^*D\cdot C.$$
where $M$ is a free $\Z$-module generated by finitely many Cartier divisors on $X$ 
and $N$ is the free $\Z$-module generated by the curves in $Y$. 
Since $Y$ has more curves than $X$, 
we can establish an appropriate boundedness result (Lemma~\ref{bounded-lemma}).

In this subsection, we often use the following notation. 

\begin{nota}\label{notation-lc}
Let $k$ be a separably closed field of characteristic $p>0$. 
Let $X$ be a projective $\Q$-factorial log canonical surface over $k$. 
%and let $f:X \to S$ be a projective morphism to a separated scheme of finite type over $k$. 
Let $B_1, \cdots, B_{\ell}$ be prime divisors. 
Fix $m\in\Z_{>0}$ such that $mK_X$ and all $mB_i$ are Cartier. 
Set 
$$M:=\Z(mK_X) \oplus \bigoplus_{i=1}^{\ell} \Z(mB_i).$$
Let 

\[
\mathcal L:=\left \{D=\sum_{i=1}^{\ell}b_iB_i \in M_{\R}\,
\middle| 
\, (X, D)\,\,{\rm is\,\,log\,\,canonical}\right \}.
\]

Let $Y$ be the normalisation of $(X\times_k \overline{k})_{\red}$ and 
let $g:Y \to X$ be the induced morphism. 
Note that $Y$ is $\Q$-factorial (Lemma~\ref{purely-basic}(3)). 
Let 
$$N:=\bigoplus \Z C_Y$$ 
where $C_Y$ runs over all the curves on $Y$. 
We obtain a $\Z$-bilinear homomorphism 
$$\varphi:M \times N \to \Z,\quad (D, C) \mapsto g^*D\cdot C.$$
\end{nota}

To apply Proposition~\ref{polytope-key}, we need the following result 
on the boundedness of extremal rays.

\begin{lem}\label{bounded-lemma}
We use Notation~\ref{notation-lc}. 
Then, there exists $\rho \in \Z_{>0}$ such that, 
for every extremal ray $R$ of $\overline{NE}(Y)$ spanned by a curve, 
there is a curve $C$ on $Y$ 
such that $R=\R_{\geq 0}[C]$ and that 
$$-g^*(K_X+B)\cdot C\leq \rho$$
for every $B \in \mathcal L$. 
\end{lem}

\begin{proof}
We obtain 
$$K_Y+\Delta_Y=g^*K_X$$
for some effective $\Z$-divisor $\Delta_Y$ by Theorem~\ref{purely}. 

We see 
$\mathcal L \subset \{\sum_{i=1}^{\ell}b_iB_i \in M_{\R}\,|\, 0\leq b_i \leq 1\}.$ 
Let $g^*B_i=\beta_iB_i'$, where $B_i'$ is the prime divisor and $\beta_i\in\Z_{>0}$. 
Set $\beta:=\max_{1\leq i\leq \ell} \beta_i$ and 
$$\mathcal L':=
\left\{\sum_{i=1}^{\ell}b_i'B_i'\,\middle|\, b_i'\in \R, 0\leq b'_i \leq \beta\right\}.$$
By \cite[Lemma~3.37]{minimal}, we can find $\rho \in \Z_{>0}$ 
such that, for every extremal ray $R$ of $\overline{NE}(Y)$ spanned by a curve, 
there exists a curve $C$ on $Y$ such that 
$$-(K_Y+\Delta_Y+B')\cdot C\leq \rho$$
for every $B' \in \mathcal L'$. 
Since $B':=g^*B$ is contained in $\mathcal L'$ for every $B \in \mathcal L$, 
we obtain 
$$-g^*(K_X+B)\cdot C=-(K_Y+\Delta_Y+B')\cdot C\leq \rho.$$
This implies the assertion. 
\end{proof}

\begin{prop}\label{polytope-main}
We use Notation~\ref{notation-lc}. 
Let $\{R_t\}_{t\in T}$ be 
the family of the extremal rays of $\overline{NE}(Y)$ spanned by a curve. 
Then the set 
$$\mathcal{N}_T:=\{B\in\mathcal{L}\,|\,
g^*(K_X+B)\cdot R_t\geq 0 \,\,for\,\, every \,\,t \in T\}$$
is a rational polytope.
\end{prop}

\begin{proof}
We fix $\rho \in \Z_{>0}$ as in Lemma~\ref{bounded-lemma}. 
By Lemma~\ref{bounded-lemma}, 
for every $t\in T$, there exists a curve $C_t$ on $Y$ 
such that $R_t=\mathbb R_{\geq 0}[C_t]$ and 
that $-g^*(K_X+B)\cdot C_t\leq \rho$ for all $B\in \mathcal L$. 
Thus, the assertion follows from Proposition~\ref{polytope-key}. 
\end{proof}

Now, we prove the abundance theorem with $\mathbb{R}$-coefficients. 

\begin{thm}\label{R-abundance}
Let $k$ be a field of characteristic $p>0$. 
Let $(X, \Delta)$ be a log canonical $k$-surface, where $\Delta$ is an $\R$-divisor. 
Let $f:X \to S$ be a projective $k$-morphism to a separated scheme $S$ of finite type over $k$. 
If $K_X+\Delta$ is $f$-nef, then $K_X+\Delta$ is $f$-semi-ample. 
\end{thm}

\begin{proof}
We may assume that $k$ is separably closed. 
By replacing $X$ with its minimal resolution, we may assume that $X$ is regular. 
In particular, $X$ is $\Q$-factorial. 

Let $\Delta=\sum_{i=1}^{\ell} \delta_iB_i$ be the irreducible decomposition, 
where each $B_i$ is a prime divisor. 
%Fix $m\in\Z_{>0}$ such that $mK_X$ and all $mB_i$ are Cartier. 
Set 
$$M:=\Z K_X \oplus \bigoplus_{i=1}^{\ell} \Z B_i$$
Note that, since $X$ is regular, $K_X$ and all $B_i$ are Cartier. 
Let 
$$\mathcal L:=\{D=\sum_{i=1}^{\ell}b_iB_i \in M_{\R}\,|\, (X, D)\,\,{\rm is\,\,log\,\,canonical}\}.$$

Let $\{R_t\}_{t\in T}$ be the set of the extremal rays of 
$\overline{NE}(X/S)$ spanned by proper $k$-curves $C$ such that $f(C)$ is one point. 
By \cite[Theorem 2.14]{MMP}, 
we obtain 
\begin{eqnarray*}
\mathcal N_T&:=&\{B\in \mathcal L \,|\,(K_X+B)\cdot R_t\geq 0\,\,{\rm for\,\,every}\,\,t\in T\}\\
&=&\{B\in \mathcal L \,|\,\,\,K_X+B\,\,{\rm is}\,\,f{\rm \mathchar`-nef}\}.
\end{eqnarray*}
If $\dim f(X)=0$, $\mathcal N_T$ is a rational polytope by Proposition~\ref{polytope-main}. 
If $\dim f(X)>0$, then the $f$-nef cone is a rational polytope by \cite[Lemma 2.13]{MMP}, 
hence also $\mathcal N_T$ is a rational polytope. 
In any case, $\mathcal N_T$ is a rational polytope. 

Since $\Delta\in \mathcal N_T$, 
we can find $\Q$-divisors $\Delta_1, \cdots, \Delta_{\ell}$  
such that $\Delta_i\in \mathcal N_T$ for all $i$ and 
that $\sum_{1 \leq i \leq \ell} r_i\Delta_i=\Delta$ 
where positive real numbers $r_i$ satisfy $\sum_{1 \leq i \leq \ell} r_i=1$. 
Thus we have 
$$K_X+\Delta=\sum_{1\leq i \leq \ell} r_i(K_X+\Delta_i)$$
and $K_X+\Delta_i$ is $f$-nef. 
By Theorem~\ref{relative-abundance}, $K_X+\Delta_i$ is $f$-semi-ample. 
\end{proof}

%% file: section5.tex
\appendix
\def\thesection{A}
\section{Shokurov polytopes in convex geometry}\label{section-convex}

In this section, we summarise some results of Shokurov polytopes 
in a setting of convex geometry. 
We fix the following notation. 

\begin{nota}\label{polytope-notation}
Let $M$ and $N$ be torsion free $\Z$-modules and let 
$$\varphi:M \times N \to \Z$$
be a $\Z$-bilinear homomorphism. 
For $D\in M$ and $C \in N$, 
we write 
$$\varphi(D, C)=D\cdot C=C\cdot D$$
by abuse of notation. 
Assume that $\dim_{\R}M_{\R}<\infty$, i.e. 
$M$ is a free $\Z$-module of finite rank. 
% ($N$:curve/not(?) nume. $M$: fin sub space of divisor space (not nume!!)). 
Fix a rational polytope $\mathcal L \subset M_{\R}$. 
%such that the origin $0\in M$ is a vertex of $\mathcal L$. 
Fix an $\R$-linear basis of $M_{\R}$ and 
we denote the sup norm with respect to this basis by $||\bullet||$. 
\end{nota}

\begin{lem}\label{polytope-lemma}
We use Notation~\ref{polytope-notation}. 
Fix $K \in M_{\Q}$, $\Delta\in\mathcal L$, and $\rho \in \R_{>0}$. 
Then, there exist positive real numbers $\epsilon, \delta>0$, 
depending on $K$, $\Delta$, and $\rho$, which satisfy the following properties. 
\begin{enumerate}
\item{For every $C \in N$ such that $-(K+B)\cdot C \leq \rho$ for all $B\in \mathcal{L}$, 
if $(K+\Delta)\cdot C>0$, then $(K+\Delta)\cdot C>\epsilon$.}
\item{If $C\in N$ and $B_0\in \mathcal{L}$ satisfy $||B_0-\Delta||<\delta$, $(K+B_0)\cdot C\leq 0$, and 
$-(K+B)\cdot C \leq \rho$ for all $B\in \mathcal{L}$, then $(K+\Delta)\cdot C\leq 0$.} 
\end{enumerate}
\end{lem}

\begin{proof}
The assertions follow from the same proof as in \cite[Proposition~3.2(1)]{Birkar1} or \cite[Theorem 4.7.2(1)]{Fuj}. 
The assertion (2) follows from the same proof as in 
\cite[Theorem 4.7.2(2)]{Fuj} 
(although the statement of 
\cite[Theorem 4.7.2(2)]{Fuj} is the same as the one of 
\cite[Proposition~3.2(2)]{Birkar1}, 
\cite[Theorem 4.7.2(2)]{Fuj} fixes minor errors appearing in 
\cite[Proposition~3.2(2)]{Birkar1}). 
\end{proof}

\begin{prop}\label{polytope-key}
We use Notation~\ref{polytope-notation}. 
Let $K \in M_{\Q}$ and $\rho \in \R_{>0}$. 
Fix a subset $\{C_t\}_{t\in T} \subset N$ 
such that 
$$-(K+B)\cdot C_t\leq \rho$$ 
for every $t \in T$ and every $B\in \mathcal{L}$. 
For any subset $T' \subset T$, we define  
$$\mathcal{N}_{T'}:=\{B\in\mathcal{L}\,|\,(K+B)\cdot C_t\geq 0 \,\,for\,\, every \,\,t \in T'\}.$$
Then there exists a finite subset $S \subset T$ such that 
$$\mathcal{N}_{T}=\mathcal N_{S}.$$
In particular $\mathcal{N}_{T}$ is a rational polytope. 
\end{prop}

\begin{proof} 
We can apply the same argument as in 
\cite[Proposition~3.2(3)]{Birkar1} by using Lemma~\ref{polytope-lemma} 
instead of \cite[Proposition~3.2(1)(2)]{Birkar1}. 
\end{proof}